\documentclass[a4paper,10pt]{article}
\usepackage{makeidx}
\usepackage[english]{babel}
\usepackage[latin1]{inputenc}
\usepackage{amsfonts}
\usepackage{bm}
\usepackage{amssymb}
\usepackage{latexsym}
\usepackage{amsmath}
\usepackage{amscd}
\usepackage{amsthm}
\usepackage{rotating}
\usepackage{graphicx}
\usepackage{pifont}
\usepackage{curves}
\usepackage{fancyhdr}
\usepackage{epsfig}
\usepackage{pstricks}
\usepackage{pst-tree}
\usepackage{subfigure}
\usepackage{epic}
\usepackage{alltt}
\usepackage[all]{xy}
\usepackage{appendix}

\newcommand{\E}{\mathbb{E}}

\newcommand{\Co}{\mathcal{C}}

\newcommand{\M}{\mathcal{M}}
\newcommand{\N}{\mathbb{N}}

\newcommand{\R}{\mathbb{R}}
\newcommand{\X}{\mathcal{X}}
\newcommand{\eg}{\textit{e.g. }}
\newcommand{\ie}{\textit{i.e.} }
\newcommand{\etal}{\textit{et al.} }
\newcommand{\wrt}{w.r.t. }

\newcommand{\lbrac}{\left[\!\left[}
\newcommand{\rbrac}{\right]\!\right]}

\def\N{\mathbb{N}}

\def\R{\mathbb{R}}
\def\E{\mathbb{E}}

\def\M{\mathcal{M}}

\def\ind{{\mathchoice {\rm 1\mskip-4mu l} {\rm 1\mskip-4mu l}
{\rm 1\mskip-4.5mu l} {\rm 1\mskip-5mu l}}}

\theoremstyle{plain}

\newtheorem{prpstn}{Proposition}

\newtheorem{crllr}{Corollary}

\theoremstyle{definition}

\newtheorem{rmrk}{Remark}

\title{Stochastic and deterministic models for age-structured populations with genetically variable traits}

\author{Régis Ferrière, Viet Chi Tran}

\begin{document}

\maketitle

\begin{abstract}Understanding how stochastic and non-linear deterministic processes interact is a major challenge in population dynamics theory. After a short review, we introduce a stochastic individual-centered particle model to describe the evolution in continuous time of a population with (continuous) age and trait structures. The individuals reproduce asexually, age, interact and die. The 'trait' is an individual heritable property (d-dimensional vector) that may influence birth and death rates and interactions between individuals, and vary by mutation. In a large population limit, the random process converges to the solution of a Gurtin-McCamy type PDE. We show that the random model has a long time behavior that differs from its deterministic limit. However, the results on the limiting PDE and large deviation techniques \textit{à la} Freidlin-Wentzell provide estimates of the extinction time and a better understanding of the long time behavior of the stochastic process. This has applications to the theory of adaptive dynamics used in evolutionary biology. We present simulations for two biological problems involving life-history trait evolution when body size is plastic and individual growth is taken into account.
\end{abstract}

This paper arose from the CANUM 2008 conference, where a mini-symposium entitled \textit{Hybrid methods} was organized by Madalina Deaconu and Tony Lelièvre. This work was supported by the French Agency for Research through the ANR MAEV and by the NSF FIBR award EF0623632.\\



\bigskip

\par In so-called structured populations, individuals differ according to variables that influence their survival and reproduction abilities. These variables are morphological, physiological or behaviorial traits that can be inherited from parents to offspring and may vary due to genetic mutation. Other types of variables include ages, which are increasing functions of time, sex, or spatial location. This paper considers the dynamics of populations with trait and age structures.\\
Age variables sum up individuals past histories and it may be biologically relevant to take several ages into account to model the present dynamics. Examples include the \textit{physical age} (the time since birth), the \textit{biological age} (the intrinsic physiological stage that can grow nonlinearly in time), the age or stage of an illness, and the time since maturation. Taking age structure into account is required to study the evolution of individuals' life histories. For example, the salmon \textit{Oncorhynchus gorbuscha} breeds when about 2 year-old, after returning to native freshwater streams, and most individuals die after this unique reproduction episode. This striking life history remains incompletely understood from an evolutionary point of viex, due in part to the lack of appropriate modelling tools. Moreover, age and trait structures interact as they change through time. This is because the age of expression of a trait and the life history shaped by the trait will influence the population dynamics, and hence the selection pressure on the trait, and its long term adaptive evolution.

\par We start with a brief review of the literature on stochastic and deterministic models of populations with age structure in Section \ref{sectionrevue}. In Section \ref{sectionrenorm1}, the links between stochastic and deterministic approaches are established: from a stochastic individual-based model (IBM), we recover the PDEs of demography \textit{via} a large population limit. Section \ref{sectionextinction} proposes some illustrations of how deterministic and stochastic approaches complement each other. The PDEs help study the measure-valued stochastic differential equations (SDEs) that arise from stochastic IBM. The SDEs in turn provide results that may not hold in large populations. An application to the study of extinction time is presented. Using these results, a second application to the so-called adaptive dynamics theory of biological evolution is developed in Section \ref{sectionAD}. This application highlights the importance of probabilistic modelling of point events. We address some interesting biological applications through simulations in Section \ref{sectionsimu}.

\section{Populations with continuous age (and trait) structure: a brief review}\label{sectionrevue}

In this brief review of the literature, we focus on population models with continuous age and time. For discrete-time models with age or stage classes, we refer the interested reader to \cite{caswell, charlesworth, thieme}.

\subsection{Deterministic models}\label{sectionrevuedet}

Population models with continuous age and time that generalize the equations of Malthus \cite{malthus} and Verhulst's well-known \textit{logistic equation} \cite{verhulst} have often been written as PDEs (\eg \cite{charlesworth, hoppensteadt, perthametransportbio, thieme, webb}).\\

\par PDE models for populations structured by a scalar age have been introduced by Sharpe and Lotka \cite{sharpelotka}, Lotka \cite{lotka}, McKendrick \cite{mckendrick} and Von Foerster \cite{vonfoerster}: $\forall t\geq 0,\, \forall a\geq 0,$
\begin{align}
 \frac{\partial n}{\partial t}(a,t)+  \frac{\partial n}{\partial a}(a,t)= -d(a)n(a,t),\quad
n(0,t) = \int_0^{+\infty} b(a)n(a,t)da,\quad n(a,0) =  n_0(a),\label{mckendrickvonfoerster}
\end{align}where $d$, $b$ and $n_0$ are nonnegative continuous functions (death and birth rates, initial age distribution) with $b$ bounded and $n_0$ integrable. Equation (\ref{mckendrickvonfoerster}) is known as McKendrick-Von Foerster PDE. It is a transport PDE that describes the aging phenomenon, with death terms and births on the boundary $a=0$. A shortcoming of these models is that they do not take into account regulations of the population size by interference interactions or environmental limitations.
\par The first nonlinear models for populations with age structure have been introduced by Gurtin and MacCamy \cite{gurtinmaccamy}. The rates $b(a)$ and $d(a)$ of (\ref{mckendrickvonfoerster}) become functions $b(a,N_t)$ and $d(a,N_t)$ of the scalar age $a\in \R_+$ and population size $N_t=\int_0^{+\infty} n(a,t)da\in \R_+$ at time $t$. Thus, the nonlinearity introduces a \textit{feed-back} term that controls the population size. A particular case involves
\begin{equation}d(a,N_t)=d(a)+\eta N_t\label{tauxmortmarcati}\end{equation}where $d(a)$ is the natural death rate at age $a$ and $\eta N_t$ is the logistic competition term. The stationary solutions, their global stabilities and estimates of the rates of decay were studied by Marcati \cite{marcati}. Stability conditions for a more general class of birth and death rates have been obtained by Farkas \cite{farkas} and Farkas \cite{farkas2}. In \cite{farkas}, functions $b(a,N_t)$ and $d(a,N_t)$ are taken in $\Co^{0,1}_b(\R_+^2,\R_+)$ with the assumption: $\exists \bar{d}>0,\, \forall (a,N)\in \R_+^2,\, d>\bar{d}$, and the characteristic equation whose roots determine the stability of the stationary solution is derived. In \cite{farkas2}, sufficient conditions for the stability of the stationary solutions are obtained for the cases where $b(a,N_t)=b_1(a)b_2(N_t)$, in which the age-dependent fertility is weighted by the interactions with other individuals.
\par The logistic competition term assumes that each competitor experiences the same interactions. Busenberg and Iannelli \cite{busenbergiannelliseparable} considered death rates that take the form $d(a)+F(t,N_t,S_1(t),\dots S_m(t))$ where $\forall i\in \lbrac 1, m\rbrac,\, S_i(t)=\int_0^{+\infty}\gamma_i(a) n(a,t)da$, $F$ is continuous in $t$ and Lipschitz continuous with respect to the other variables, and the $\gamma_i$ are positive continuous bounded functions that describe the competitive pressure exerted by an individual depending on its age. However the interaction term remains non-local (and thus finite-dimensional). Local interactions mean that the competition terms depend on integrals of the form $\int_0^{+\infty} U(a,\alpha) n(\alpha,t)d\alpha$ where $U(a,\alpha)$ is the interaction exerted by an individual of age $\alpha$ on an individual of age $a$. This case is covered by the work of Webb \cite{webb}, who proposes PDEs with coefficients of the form $b(a,n(.,t))$ and $d(a,n(.,t))$, where $n(.,t)$ belongs to the Banach space $L^1(\R_+,\R_+)$ of integrable functions \wrt $a\in \R_+$. 
\par More recently, entropy methods have extended the notion of relative entropy to equations that are not conservation laws. This approach (Michel \etal \cite{michelmischlerperthame}, Mischler \etal \cite{mischlerperthameryzhik}, Perthame \cite{perthametransportbio}, Perthame and Ryzhik \cite{perthameryzhik}) provided new proofs of existence of solutions, and new insights into long time behavior for instance.
\par There are many examples in the literature where additional structures other than age-based have been taken into account. The population can be divided into a finite number of classes: see Webb \cite{webb} for a survey on multi-type models, including prey-predator and epidemiological models (\eg \cite{busenbergiannellithieme, iannellimilnerpugliese, lafayelanglais}). In epidemics models, considering the classes of susceptible, infectious and/or removed individuals with different interactions between these groups, leads to different propagation dynamics of the infection. Age-structure is important as infection or detection rate may depend for instance on the time since the first infection.
\par To our knowledge, populations with age structure and spatial diffusions (where motion is described with a second order differential operator) were first studied by Langlais \cite{langlais,langlais2,langlais3} and Busenberg and Iannelli \cite{busenbergiannelli}. Structure variables with values in a non-finite state-space other than position have been considered by Rotenberg \cite{rotenberg} for instance  (for multitype populations, see \eg \cite{webb}). In \cite{rotenberg}, the maturation velocity $x\in \X=]0,1]$ is defined as a trait that can change at birth and during the life of an individual. For an individual born at time $c\in \R$ and with traits $x_i\in ]0,1]$ on $[t_i,t_{i+1}[$ ($i\in \N$, $t_0=c$ and $(t_i)_{i\in \N}$ are times of trait change in $[c,+\infty[$), the biological age at $t>c$ is $a=\sum_{i\in \N}x_i(t_{i+1}\wedge t-t_i\wedge t)$. Here, the biological age $a$ differs from the physical age $t-c$ as soon as $ \exists i\in \N,\, t_i< t $ and $x_i\not= 1$. A nonlinear version of this model has been studied by Mischler \etal \cite{mischlerperthameryzhik}.

\subsection{Stochastic models}\label{reviewsto}

\par Age structured branching processes that generalize the Galton-Watson process \cite{watsongalton} have been studied by Bellman and Harris \cite{bellmanharris, harrislivre}, and then Athreya and Ney (\cite{athreyaney} Chapter IV). In these non Markovian models, the lifespan of an individual does not follow an exponential law. Upon death, a particle is replaced by a random number of daughter particles, with a reproduction law that does not depend on the age of the mother, nor on the state of the population.
\par The assumptions of births at the parent's death and of independence between the reproduction law and the age of the parent are biologically restrictive. Kendall \cite{kendall}, Crump and Mode \cite{crumpmode, crumpmode2}, Jagers \cite{jagers,jagerslivre}, Doney \cite{doney}, study birth and death processes in which an individual can give birth at several random times during its life, with rates that may depend on its age. These processes are called age-structured birth and death processes. More recently Lambert \cite{lambert} used contour processes to study the properties of splitting trees which are formed by individuals with independent and identically distributed (\textit{i.i.d.}) lifespans and who give birth at the same constant rate during their lives. Lambert's results provide a new interpretation of the link between branching processes and Lévy processes \cite{legalllejan,popovic}.
\par In the two types of models introduced above, individuals alive at the same time are independent, which is also a biologically restrictive assumption. Wang \cite{wangbio}, Solomon \cite{solomon} have considered birth and death processes in which the lifespans of individuals are independent, but in which the birth rate and reproduction law of each individual depend on the state of the population. Oelschläger \cite{oelschlager}, Jagers \cite{jagers}, Jagers and Klebaner \cite{jagersklebaner} generalize these models by including interactions in both the birth and death rates, but these rates remain bounded. In the models mentioned in this paragraph, the population at time $t$ is discrete and represented by a point measure of the form $Z_t(d\alpha)=\sum_{i=1}^{N_t}\delta_{a_i(t)}(d\alpha)$. Each individual is described by a Dirac mass charging its age $a_i(t)$, and $N_t$ is the size of the population. The birth and death rates considered by \cite{jagers, jagersklebaner, oelschlager} have the form $b(a,\langle Z_t,U(a,.)\rangle)$ and $d(a,\langle Z_t, V(a,.)\rangle)$, where $U$ and $V$ are bounded interaction kernels, $b$ and $d$ are bounded functions, and we use the notation $\langle Z,f\rangle = \int_0^{+\infty}f(\alpha)Z_t(d\alpha)$. In \cite{jagersklebaner,jagers} for instance, the interactions vanish in the limit: $\lim_{\langle Z,1\rangle \rightarrow +\infty}b(a,\int_0^{+\infty}U(a,\alpha)Z(d\alpha))=b(a)$ and the same is true for the death term, so that the authors recover population behaviors that are independent of population size in the limit of large populations. These models thus exclude logistic interactions. Such interactions are included in the models considered in \cite{trangdesdev}.
\par Interactions between individuals imply that the key \textit{branching property} ceased to hold, and classical approaches based on generating or Laplace functions (\eg \cite{athreyaney}) do not apply anymore. However, in a large population limit, a law of large numbers can be established \cite{wangbio, solomon, oelschlager, trangdesdev} whereby a limiting deterministic process identified as a weak solution of a generalization of the McKendrick-Von Foerster PDE (\ref{mckendrickvonfoerster}) is obtained. We present a synthesis of these results in Section \ref{sectionrenorm1}. Once the link between PDEs and stochastic processes is established, results for PDEs prove to have interesting implications for the stochastic process, even if the behaviors of the stochastic model and of its corresponding deterministic limit differ (see Sections \ref{sectionextinction} and \ref{sectionAD}).
\par Let us mention that there are other large-population scalings that have been the focus of abundant research. Superprocess limits of age-structured branching processes (where the mother gives birth at her death) are obtained by rescaling the lifespan and mass of the particles, by increasing their number and by modifications of the birth rate or reproduction law (\eg Dynkin \cite{dynkin91}, Kaj et Sagitov \cite{kajsagitov}). For age-structured birth and death processes (where the mother gives birth at random times during her life), similar approaches were taken by Dawson \textit{et al.} \cite{dawsongorostizali} and Bose and Kaj \cite{bosekaj, bosekaj2}.
\par Multitype populations with age structure have been considered by many authors (\eg Athreya and Ney \cite{athreyaney}, Chap.V.10). A stochastic model of epidemics with age structure has been considered for instance in Clémençon \etal \cite{arazozaclemencontran}. Stochastic models of populations with age and trait structures have been addressed by Jagers \cite{jagersgeneralbranching,jagerscoupling} and Méléard and Tran \cite{meleardtran} for instance. In \cite{meleardtran}, the mutation rate that is responsible for generating variation in the trait is
decreased so that the mutations become more and more rare. As a consequence, the timescales of mutations and demography (births and deaths) become separate, which is a typical assumption of the evolutionary theory of adaptive dynamics (see Section \ref{sectionAD}).

\section{From stochastic individual-based models of age-structured populations to deterministic PDEs}\label{sectionrenorm1}

\subsection{A microscopic model}

We study a microscopic stochastic model, for which the dynamics is specified at the level of individuals. The model takes into account trait and age dependence of birth and death rates, together with interactions between individuals. This process generalizes the approach of Fournier and Méléard \cite{fourniermeleard}, Champagnat \textit{et al.} \cite{champagnatferrieremeleard, champagnatferrieremeleard2}. When the population is large, we establish a macroscopic approximation of the process, which describes the evolution at the scale of the population (individual paths are lost). These results are taken from \cite{trangdesdev} and allow to link the probabilistic point of view to the deterministic approach reviewed in Section 1.
\par The IBM represents the discrete population at time $t\geq 0$ by a point measure \begin{equation}
Z_t=\sum_{i=1}^{N_t}\delta_{(x_i(t),a_i(t))},\label{Zintro}
\end{equation}where each individual is described by a Dirac mass $(x,a)\in \widetilde{\X}:=\mathcal{X}\times \R_+$, $\X\subset \R^{d}$ being the trait space and $\R_+$ the age space. We denote by $N_t$ the number of individuals alive at time $t$. The space of point measures on $\widetilde{\X}$ is denoted by $\mathcal{M}_P(\widetilde{\X})$ and embedded with the topology of weak convergence. For a function $f$, we denote by $\langle Z_t,f\rangle$ the integral $\int_{\widetilde{\X}}f(x,a)Z_t(dx,da)=\sum_{i=1}^{N_t}f(x_i(t),a_i(t))$. An individual of trait $x$ and age $a$ in a population $Z\in \mathcal{M}_P(\widetilde{\X})$ reproduces asexually, ages and dies. When reproduction occurs, the trait is transmitted to offspring unless a mutation occurs. The mechanism is modelled as follows:
\begin{itemize}
\item Birth rate is $b(x,a)\in \R_+$. With probability $p\in [0,1]$, the new individual is a mutant with trait $x+h$ where $h$ is chosen in the probability distribution $k(x,a,h)\,dh$ (dependent upon the characteristics of the parent) where $dh$ is the Lebesgue measure,
\item Death rate is $d(x,a,ZU(x,a))\in \R_+$. It depends on the trait and age of the individual but also on the measure $Z$ describing the whole population. The function $U\,:\,\widetilde{\X}^2\mapsto \R^m$ is an interaction kernel: each of the $m$ component of $U((x,a),(y,\alpha))$ describes a different interaction of $(y,\alpha)$ on $(x,a)$, and $ZU(x,a)=\int_{\widetilde{\X}}U((x,a),(y,\alpha))Z(dy,d\alpha)$.
\item Here, aging velocity is 1.\end{itemize}

\noindent\textbf{Assumption (H1)} The birth rate $b(x,a)$ is assumed to be continuous and bounded by a positive constant $\bar{b}$. The function $k(x,a,h)$ is assumed to be bounded by a positive constant $\bar{k}$. The death rate is assumed to be continuous in $(x,a)$, Lipschitz continuous in the interaction term, lower bounded by a strictly positive constant and upper bounded by $\bar{d}(1+\langle Z,1\rangle)$ with $\bar{d}>0$.

\begin{rmrk}1. More general models are considered in \cite{chithese, trangdesdev}, where vectorial ages with nonlinear aging velocities are allowed.\\
2. For appropriate choices of $b$ and $d$, some of the models described in Section \ref{reviewsto} can be recovered. If $b(x,a)=b(a)$ and $d(x,a)=d(a)$, we recover the age-dependent birth and death processes of \cite{crumpmode, crumpmode2, doney, jagers, kendall}. The choice $b(x,a)=b(a)$ and $d(x,a,ZU(x,a))=d(a)+\eta \langle Z,1\rangle$ gives the rates used by \cite{gurtinmaccamy}. Considering bounded functions $d$ and vectorial functions $U$ leads to models similar to those considered in \cite{jagersklebaner, oelschlager}.
\end{rmrk}

\par Following \cite{fourniermeleard,champagnatferrieremeleard, champagnatferrieremeleard2, trangdesdev}, we describe the evolution of $(Z_t)_{t\in \R_+}$ by a SDE driven by a Poisson point process. The mathematical difficulty lies in the fact that the birth and death rates depend on time through age.
\par Let $Z_0\in \mathcal{M}_P(\widetilde{\X})$ be a random variable (r.v.) such that $\mathbb{E}\left(\langle Z_0,1\rangle \right)<+\infty,$ and let $Q(ds,di,d\theta,dh)$ be a Poisson point measure (P.P.M.) on $\R_+\times \mathcal{E}:=\R_+\times \N^*\times \R_+\times \X$ with intensity $q(ds,di,d\theta,dh):=ds\otimes n(di)\otimes d\theta\otimes \,dh$ and independent of $Z_0$ (for the general theory of P.P.Ms, see for instance \cite{ikedawatanabe} Chapter I).
\par Let us denote by $X_i(t)$ and $A_i(t)$ the trait and age of the $i^{\mbox{th}}$ individual at time $t$, the individuals being ranked in the lexicographical order on $\R^{d}\times \R_+$ (see \cite{fourniermeleard,trangdesdev} for a rigorous notation).
\begin{align}
Z_t   =  & \sum_{i=1}^{N_0} \delta_{(X_i(0),A_{i}(0)+t)}
+
 \int_{0}^t \int_{\mathcal{E}} \mathbf{1}_{\{i\leq N_{s_-}\}}\left[ \delta_{(X_i(s_-),t-s)} \mathbf{1}_{\{0\leq
\theta< m_1(s,Z_{s_-},i,h)\}}\right.\nonumber\\
 + &     \delta_{(X_i(s_-)+h,t-s)} \mathbf{1}_{\{m_1(s,Z_{s_-},i,h)\leq \theta<
m_2(s,Z_{s_-},i,h)\}}\nonumber\\
- & \left.  \delta_{(X_i(s_-),A_{i}(s_-)+t-s)} \mathbf{1}_{\{m_2(s,Z_{s_-},i,h)\leq \theta <
m_3(s,Z_{s_-},i,h)\}} \right] Q(ds, di,d\theta,dh),\label{EDS}
\end{align}where:
\begin{eqnarray*}
m_1(s,Z_{s_-},i,h) & = & (1-p)b(X_i(s_-),A_i(s_-))k(X_i(s_-),A_i(s_-),h)\\
m_2(s,Z_{s_-},i,h) & = & m_1(s,Z_{s_-},i,h) +   p\,b(X_i(s_-),A_i(s_-))k(X_i(s_-),A_i(s_-),h)\\
m_3(s,Z_{s_-},i,h) & = & m_2(s,Z_{s_-},i,h) +  d(X_i(s_-),A_i(s_-),Z_{s_-}U(X_i(s_-),A_i(s_-)))k(X_i(s_-),A_i(s_-),h).
\end{eqnarray*}
The interpretation is as follows. To describe the population at any given time $t$, we start with the individuals present at $t=0$. Then we add all the ones that were born between 0 and $t$ and finally we delete the Dirac masses corresponding to individuals who died between 0 and $t$. The individuals of the initial condition have at time $t$ an age increased by $t$ (first term of (\ref{EDS})). When a birth takes place at time $s$, in the measure describing the population at time $t$ we add a Dirac mass at age $t-s$ (the age expected for this newborn individual, second and third term of (\ref{EDS})). When an individual of age $a$ dies at time $s$, in the measure describing the population at time $t$ we suppress a Dirac mass at age $a+t-s$ (the age this individual would have had if it had survived). As the rates vary with time, we use a P.P.M. with an intensity that upper-bounds the age-dependent rates and use an acceptance-rejection procedure (the indicators in $\theta$) to recover the rates that we need.
\par We emphasize that exact simulation of this SDE can be done easily (see Section \ref{sectionsimu}). The algorithm formally accounts for the computer programs that many biologists use to run simulations of individual-based population models (\eg \cite{deangelisgross}).

\begin{prpstn}\label{propexistence}
Under \textbf{(H1)}, for every given Poisson point measure $Q$ on $\R_+\times \mathcal{E}$ with intensity measure $q$, and every initial condition $Z_0$ such that $\mathbb{E}\left(\langle Z_0,1\rangle\right)<+\infty$, there exists a unique strong solution to SDE (\ref{EDS}) in $\mathbb{D}(\R_+,\mathcal{M}_P(\widetilde{\X}))$. The solution is a Markov process with infinitesimal generator given for all $F_f(Z)=F(\langle Z,f\rangle)$ with $f\in \Co^{0,1}(\widetilde{\X},\R)$, $F\in \Co^1(\R,\R)$, and $Z\in \mathcal{M}_P(\widetilde{\X})$ by
\begin{align}L F_f (Z) = \int_{\widetilde{\mathcal{X}}} & \left[
 \partial_{a} f(x,a)F'\left(\langle Z, f\rangle \right)
+\left(F_f\left(Z+\delta_{(x,0)}\right)-F_f\left(Z\right)\right)
b(x,a)(1-p)\right. \nonumber\\
 + &
\int_{\R^d}\left(F_f\left(Z+\delta_{(x+h,0)}\right)-F_f\left(Z\right)\right)
b(x,a)p\,k(x,a,h)dh\nonumber\\
 +  & \left. \left(F_f\left(Z-\delta_{(x,a)}\right)-F_f\left(Z\right)\right)d(x,a,Z U(x,a))\right]
Z(dx,da).
\end{align}
\end{prpstn}

For the proof, see Propositions 2.2.5, 2.2.6 and Theorem 2.2.8 in \cite{chithese}. The result is based on moment estimates and uses the algorithmic construction of (\ref{EDS}) that will be detailed in Section \ref{sectionsimu}. The behaviour of the solution is far more difficult to study. For this reason, and also in order to link the SDE (\ref{EDS}) with the PDEs that are classically introduced in demography (see Section \ref{sectionrevuedet}), we consider a large population limit. Thus we let the size of the initial population grow to infinity proportionally to an integer parameter $n$ while individuals and interaction intensities are assigned a weight $1/n$. Under proper assumptions on the initial condition, the sequence of rescaled processes $(Z^n)_{n\in \N}$ converges, when $n\rightarrow +\infty$, to the solution of a deterministic equation identified as the weak form of a PDE ; this PDE generalizes the deterministic models described in Section \ref{sectionrevuedet}.

\subsection{Large population limit}\label{sectionrenorm}

The large population renormalization that we consider is inspired by the work of Fournier and Méléard \cite{fourniermeleard}. Full details and proofs can be found in \cite{trangdesdev,chithese}. We let the size of the initial population tend to infinity proportionally to the parameter $n\in \N^*$ and renormalize the weights of the individuals and their interaction by $1/n$. Following on Metz \etal \cite{metzgeritzmeszenajacobsheerwaarden}, $n$ was named "system size" by Champagnat \etal \cite{champagnatferrieremeleard}. More precisely, we consider a sequence $(Z^n_t,\,t\in \R_+)_{n\in \N^*}$ such that $nZ^n$ satisfies SDE (\ref{EDS}) with $U^n=U/n$ so that the interaction term $ZU(x,a)=\sum_{i=1}^{N_t}U((x,a),(x_i,a_i))$ is replaced with $\frac{1}{n}\sum_{i=1}^{N^n_t}U((x,a),(x_i,a_i))=Z^n U(x,a)$ and such that:\\

\noindent \textbf{Assumption (H2)}: The sequence $(Z^n_0)_{n\in \N^*}=\big(\frac{1}{n}\sum_{i=1}^{N^n_t}\delta_{(X_i(0),A_i(0))}\big)$ converges in probability to the measure $\xi_0$ belonging to the space $\mathcal{M}_F(\widetilde{\X})$ of finite measures embedded with the weak convergence topology.\\

This large population rescaling can be understood as reflecting resource limitation as the size of the population increases: we have to assume that the biomass of individuals decreases if the system is to remain viable. As the individuals become smaller, their interactions decrease proportionally to their mass.

\begin{prpstn}\label{proplgn}Under \textbf{(H1)} and \textbf{(H2)},\\
(i) If $\sup_{n\in \N^*}\E\big(\langle Z^n_0,1\rangle^{2}\big) <+\infty$, then for every $f\in \Co^{1,1,0}(\R_+\times \widetilde{\X},\R)$ the process
\begin{align}
M_t^{n,f}  = &    \langle Z^n_t, f(.,.,t) \rangle - \langle Z^n_0, f(.,.,0) \rangle - \int_0^t \int_{\widetilde{\X}}\left[\partial_{a} f(x,a,s)+ \partial_sf(x,a,s)+ f(x,0,s)b(x,a)(1-p) \right.\nonumber\\
+  &  \left. \int_{\R^d} f(x+h,0,s)b(x,a)p\,k(x,a,h)dh
-  f(x,a,s)d(x,a,Z^n_sU(x,a))\right]Z_s^n(dx,da)\, ds,\label{pbm}
\end{align}
is a square integrable martingale with the following quadratic variation process:
\begin{align}
\langle M^{n,f} \rangle_t  =
\frac{1}{n} \int_0^t \int_{\widetilde{\X}} & \left[f^2(x,0,s)b(x,a)(1-p)+  \int_{\R^d} f^2(x+h,0,s)b(x,a)p\,k(x,a,h)dh \right.\nonumber\\
 + & \left. f^2(x,a,s)d(x,a,Z^n_sU(x,a))\right] Z_s^n(dx,da)\, ds.
\label{crochetmartingale}\end{align}
(ii) If $\exists \eta>0,\, \sup_{n\in \N^*}\E\big(\langle Z^n_0,1\rangle^{2+\eta}\big) <+\infty$, then
the sequence $(Z^n)_{n\in \N^*}$ converges in distribution in $\mathbb{D}(\R_+,\mathcal{M}_F(\widetilde{\X}))$ to the unique solution $\xi\in \Co(\R_+,\mathcal{M}_F(\widetilde{\X}))$ of the following equation: $\forall f\,:\,(s,x,a)\mapsto f_s(x,a) \in \Co^{1,0,1}(\R_+\times \widetilde{\X},\R),$
\begin{align}
\langle \xi_t, f(.,.,t) \rangle   =  & \langle \xi_0, f(.,.,0) \rangle + \int_0^t \int_{\widetilde{\mathcal{X}}}\left[ \partial_{a} f(x,a,s)+\partial_s f(x,a,s)+f(x,0,s)b(x,a)(1-p)\right.\nonumber\\
+ & \left.  \int_{\R^d}f(x+h,0,s)b(x,a)p\, k(x,a,h)dh
- f(x,a,s)d(x,a,\xi_sU(x,a))\right]\xi_s(dx,da)\,ds,\label{formefaible}
\end{align}
\end{prpstn}

This result is proved by a tightness-uniqueness argument (see \cite{fourniermeleard, trangdesdev}). Point (i) is obtained by stochastic calculus for jump processes. Heuristically, since the quadratic variation is of the order of $1/n$, it vanishes as $n\rightarrow +\infty$. Moreover, when births or deaths occur, we add or delete individuals of weight $1/n$. This explains why in the limit we obtain a continuous deterministic process.\\

\par The link between (\ref{formefaible}) and PDEs is specified in the next proposition (for proofs see Prop.3.6 and 3.7 in \cite{trangdesdev}).\\

\noindent\textbf{Assumption (H3):} $\xi_0$ admits a density $n_0(x,a)$ with respect to $dx\otimes da$ on $\widetilde{\X}$.\\

\begin{prpstn}\label{propdensite}Under \textbf{(H1)}, \textbf{(H2)} and \textbf{(H3)}, the measures $\xi_t\in \mathcal{M}_F(\widetilde{\X})$, for $t\in \R_+$, admit densities $n(x,a,t)$ with respect to $dx\otimes da$ on $\widetilde{\X}$. The family of these densities $(n(.,.,t))_{t\in \R_+}$ is a weak solution of:
\begin{eqnarray}
\partial_t n(x,a,t) & = & - \partial_{a}n(x,a,t)-
 d\left(x,a,\int_{\widetilde{\X}} U((x,a),(y,\alpha)) n(y,\alpha,t)dy\, d\alpha\right) n(x,a,t) \label{edpintropop1}\\
n(x,0,t) & = & \int_{\R_+} n(x,a,t) b(x,a)(1-p) da+
 \int_{\R^d \times \R_+}
b(x-h,a)p\,k(x-h,a,h)n(x-h,a,t)dh\,da\nonumber\\
n(x,a,0) & = & n_0(x,a).\nonumber
\end{eqnarray}
\end{prpstn}

These equations generalize the McKendrick-Von Foerster's PDEs. They describe ecological dynamics at the scale of the population (individual trajectories are lost). The densities $n(x,a,t)$, when they exist, correspond to the \textit{number density} in the sense of Desvillettes \etal \cite{desvillettesferriereprevost}, Champagnat \etal \cite{champagnatferrieremeleard2}, which describes the trait and age distributions of a "continuum" of individuals. Here, existence and uniqueness of a weak solution of (\ref{edpintropop1}) are obtained by probabilistic proofs. An interesting property of equation (\ref{formefaible}) is that function solutions do not always exist (see \cite{trangdesdev}).

\subsection{Central limit theorem}

\par In large populations, the microscopic process $Z^n$ can be approximated by the solution of PDE (\ref{edpintropop1}) that generalizes classical PDEs of demography for populations in continuous time and structured by a scalar age. In order to construct confidence intervals or to assess the quality of the approximation, it is useful to study the fluctuation process,
\begin{align}
\forall t\in [0,T],\, \forall n\in \N^*,\quad \eta_t^n(dx,da)= & \sqrt{n}\left(Z^n_t(dx,da)-\xi_t(dx,da)\right).
\end{align}
For $n\in \N^*$, $\eta^n$ is a process of $\mathbb{D}([0,T], \mathcal{M}_S(\widetilde{\X}))$. Since the space $\mathcal{M}_S(\widetilde{\X})$ of signed measures on $\widetilde{\X}$ cannot be meterized when embedded with the topology of weak convergence, we follow the works of Métivier \cite{metivierIHP} and Méléard \cite{meleardfluctuation}, and consider $\eta^n$ as a distribution-valued process.
\par Let us denote by $\Co_K^\infty$ the space of smooth functions with compact support. For a multi-index $k\in \N^{d+1}$, we denote by $D^kf$ the derivative $\partial^{|k|}f/\partial x_1^{k_1}\cdots \partial x_d^{k_d}\partial a^{k_{d+1}}$. The Banach space $C^{\beta,\gamma}$ is the space of functions $f$ of class $\Co^\beta$ such that for $k$ such that $|k|\leq \beta$, $\lim D^k f(x,a)/(1+a^\gamma)=0$ when $a\rightarrow +\infty$ or when $x\rightarrow \partial \X$, embedded with the norm $\|f\|_{C^{\beta,\gamma}}=\sum_{|k|\leq \beta}\sup_{(x,a)\in \widetilde{\X}} |D^k f(x,a)|/(1+a^\gamma)$. The Hilbert space $W_0^{\beta,\gamma}$ is the closure of $\Co_K^\infty$ with respect to the norm $\|f\|_{W_0^{\beta,\gamma}}^2=\int_{\widetilde{\X}} \sum_{|k|\leq \beta}|D^k f(x,a)|^2/(1+a^{2\gamma}) \, dx\,da$, $\Co_K^\infty$ being the space of smooth functions with compact support. We denote by $C^{-\beta,\gamma}$ and $W_0^{-\beta,\gamma}$ their dual spaces. We are going to consider the following embeddings. Let $D=[(d+1)/2]+1$,
\begin{align}
C^{3D+1,0}\hookrightarrow W_0^{3D+1,D} \hookrightarrow_{H.S.} W_0^{2D+1,2D}\hookrightarrow C^{D+1,2D}\hookrightarrow C^{D,2D} \hookrightarrow W_0^{D,3D}\hookrightarrow C^{0,3D+1}
\end{align}the second embedding being Hilbert-Schmidt (\eg \cite{adams} p.173).
\par The tightness is proved in the space $W_0^{-3D+1,D}$ using a tightness criterion for Hilbert-valued processes that is stated in \cite{meleardfluctuation} (Lemma C), and which requires the Hilbert-Schmidt embedding $W_0^{3D+1,D} \hookrightarrow_{H.S.} W_0^{2D+1,2D}$. $C^{-(D+1),2D}$ and $C^{-D,2D}$ are the spaces in which the estimates of the norm of the finite variation part of the process are established (see Section 4.4.2 of \cite{chithese}). We need these two spaces since $\partial_a$ maps $C^{D+1,2D}$ into $C^{D,2D}$. The norm of the martingale part is controlled in $W_0^{-D,3D}$. Finally, uniqueness of the limiting value is proved thanks to the embedding $ W_0^{-3D+1,D}\hookrightarrow  C^{-3D+1,0}$.\\

\noindent\textbf{Assumption (H4)}: We assume that Assumptions \textbf{(H1)}, \textbf{(H2)} and \textbf{(H3)} are satisfied and that:
 \begin{itemize}
 \item The sequence $(\eta^n_0)_{n\in \N^*}$ converges in $W_0^{-D,3D}$ to $\eta_0$ and $\sup_{n\in \N^*}\mathbb{E}\big(\|\eta^n_0\|^2_{W_0^{-D,3D}}\big)<+\infty$.
 \item $\sup_{n\in \N^*}\mathbb{E}\big(\big(\int_{\widetilde{\X}}|a|^{6D}Z^n_0(dx,da)\big)^2\big)<+\infty$.
 \item The functions $(x,a)\mapsto b(x,a)$, $(x,a)\mapsto k(x,a,x')$, $(x,a)\mapsto U((x',\alpha),(x,a))$ and $(x,a)\mapsto U((x,a),(x',\alpha))$ belong to $C^{3D+1,0}$ for almost every $(x',\alpha)\in \widetilde{\X}$.
 \item The death rate $(x,a,u)\mapsto d(x,a,u)$ is such that for all $u$, it belongs to $C^{3D+2,0}$ with a norm bounded by a polynomial in $u$. We assume that its derivative with respect to $u$ is Lipschitz continuous.
\end{itemize}

\begin{prpstn}\label{proptcl}
Let $T>0$. Under Assumptions \textbf{(H4)}, the sequence $(\eta^n)_{n\in \N^*}$ converges in distribution in $\mathbb{D}([0,T],W_0^{-3D+1,D})$ to the unique continuous solution of: $\forall t\in [0,T],\, \forall f\in W_0^{3D+1,D},$
\begin{align}
\langle \eta_t,f\rangle = & \langle \eta_0,f\rangle + W_t(f)+  \int_0^t \int_{\widetilde{\X}}\left[\partial_a f(x,a)+b(x,a)\big((1-p)f(x,0)+p\int_{\X} f(x+h,0)k(x,a,h)dh\big)\right.\nonumber\\
- & \left.d(x,a,\xi_s U(x,a))f(x,a)-\int_{\widetilde{\X}}f(y,\alpha)d_u(y,\alpha,\xi_sU(y,\alpha))U((y,\alpha),(x,a))\xi_s(dy,d\alpha)\right] \eta_s(dx,da)\, ds
\end{align}where $(W_t(f))_{t\in \R_+}$ is a continuous centered square-integrable Gaussian process with quadratic variation:
\begin{align}
\langle W(f)\rangle_t=  & \int_0^t \int_{\widetilde{\X}} \left[b(x,a)\big((1-p) f^2(x,a)+p \int_{\X}f^2(x+h,0)k(x,a,h)dh\big)\right.\nonumber\\
+  & \left. f^2(x,a)d(x,a,\xi_sU(x,a))\right]\xi_s(dx,da)\, ds.
\end{align}
\end{prpstn}

Notice that a benefit from IBMs is that the stochastic approach with point measures suits the formalism of statistical methods which involve sets of individual data (\eg \cite{arazozaclemencontran, blumtran}). This makes it possible to use the battery of statistical methods which deal with problems that may not always be treated with deterministic methods (such as missing or noisy data). The convergence and fluctuations of Propositions \ref{proplgn} and \ref{proptcl} would then provide consistence and asymptotic normality for the estimators. Calibrating the parameters of a PDE thanks to its underlying microscopic interpretation provides an alternative to deterministic approaches (\eg \cite{banksburnscliff, burnscliffdoughty}).

\subsection{Extinction in logistic age structured populations with constant trait}\label{sectionextinction}

Extinction is one of the recurrent and important issue when studying the ecology of a population. Here, the conclusions given by the deterministic and the stochastic models differ, but we will see that they still provide complementary information.
\par For the sake of simplicity, we consider here a logistic age-structured population without trait variation. We will write for instance $n_0(a)$ instead of $n_0(x,a)$ of Assumption \textbf{(H2)}. An individual of age $a$ in a population of size $N$ gives birth with rate $b(a)$ and dies with the rate $d(a)+\eta N$ as in (\ref{tauxmortmarcati}). From Proposition \ref{proplgn}, under the large population renormalization of Section \ref{sectionrenorm}, the sequence $(Z^n)_{n\in \N^*}$ converges to the unique weak solution of the Gurtin McCamy PDE with the death rate (\ref{tauxmortmarcati}). Let us recall some known facts about this PDE. We denote by $\Pi(a_1,a_2)=\exp\big(\int_{a_1}^{a_2} d(\alpha)d\alpha\big)$ the probability of survival from age $a_1$ to age $a_2$ in absence of competition.

\begin{prpstn}\label{propedpclassique}Under \textbf{(H1)} and \textbf{(H3)}, there exists a unique classical solution of class $\Co^1$ to the Gurtin-McCamy equation with death rate (\ref{tauxmortmarcati}). It is given by: $\forall a\in \R_+,\,\forall t\in \R_+$
\begin{align}
n(a,t)=\frac{N_0 v(a,t)}{1+N_0\int_0^t\int_0^{+\infty}v(\alpha,s)d\alpha\,ds}\quad \mbox{ with }v(a,t)=\left\{\begin{array}{ccc}
n_0(a-t)\Pi(a-t,a)/N_0 & \mbox{ if } & a\geq t\\
B(t-a)\Pi(a,0)& \mbox{ if } & a< t\\
\end{array}\right.\label{nat}
\end{align}where $n_0$ is the density of the initial condition, $N_0$ its $L^1$-norm and:
\begin{align*}
B(t)=B_0*\big(\sum_{n=0}^{+\infty}g^{*n}\big)(t)\mbox{ with }g(a)=b(a)\Pi(0,a)\ind_{a\geq 0} \mbox{ and } B_0(t)=\ind_{t\geq 0}\int_0^{+\infty}b(a+t)\frac{n_0(a)}{N_0}\Pi(a,a+t)da.\end{align*}
\end{prpstn}
The quantity $B(t)$ which appears in the computations can be interpreted heuristically as the number of births at time $t$. For a proof, see Webb \cite{webb} (Section 5.4).
\par Since the family $(n(a,t)da)_{t\in \R_+}$ defines a weak solution of the Gurtin-McCamy PDE in this case, we obtain the following corollary of the uniqueness property stated in Proposition \ref{proplgn} and of Proposition \ref{propedpclassique}:
\begin{crllr}
In the case of logistic age-structured population, the limiting process $\xi$ of the sequence $(Z^n)_{n\in \N^*}$ satisfies:
$\forall t\in \R_+,\, \xi_t(da)=n(a,t)da$ where $n(a,t)$ is explicitly given in (\ref{nat}).\end{crllr}We have used here the results known for the PDE to obtain an explicit expression and regularities of the density of $\xi$ obtained by probabilistic proofs (Propositions \ref{proplgn} and \ref{propdensite}).

\par As explained in Section \ref{sectionrevue}, the long time behavior of the Gurtin-McCamy PDE is well known. In particular:

\begin{prpstn}(see \cite{marcati,webb})
Under the assumption $R_0:=\int_0^{+\infty}b(a)\Pi(0,a)da >1$, which expresses the renewal of generations in absence of competition, there exists a unique nontrivial asymptotically exponentially stable steady state given by:
\begin{align}\label{nchap}
\widehat{n}(a) = \frac{\lambda_1 e^{\lambda_1 a}\Pi(0,a)}{\eta \int_0^{+\infty} e^{-\lambda_1 \alpha}\Pi(0,\alpha)d\alpha}\quad \mbox{where }\lambda_1\mbox{ satisfies}\quad
1=\int_0^{+\infty} e^{-\lambda_1 a}b(a)\Pi(0,a)da.\end{align}
\end{prpstn}

\begin{figure}[!ht]
\begin{center}
\begin{tabular}[!ht]{ccc}
\hspace{0.7cm}
\includegraphics[width=0.25\textwidth,height=0.15\textheight,angle=270,trim= 2cm 4cm 2cm 2cm]{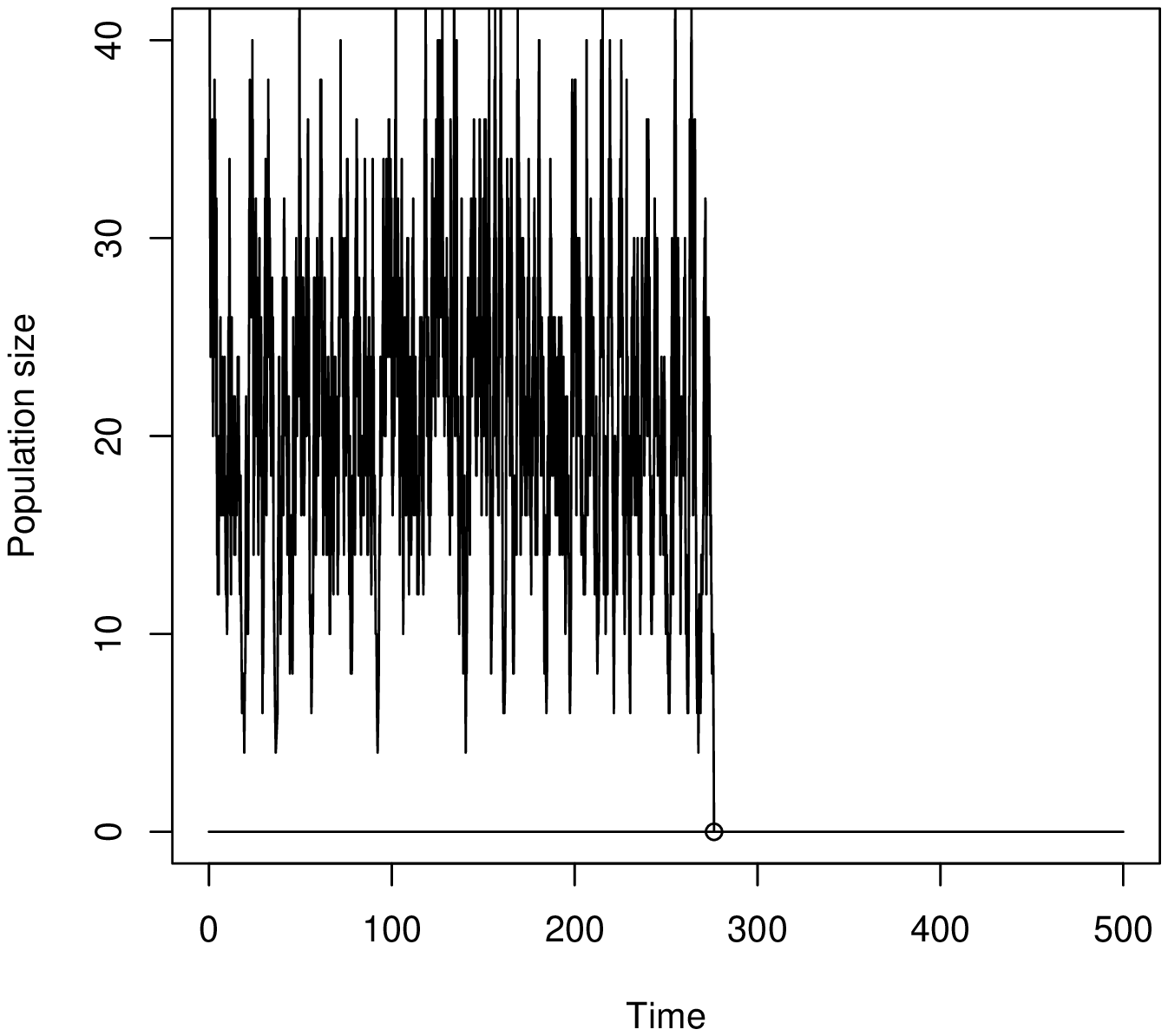} \vspace{0.2cm}&
\hspace{1.5cm}
\includegraphics[width=0.25\textwidth,height=0.15\textheight,angle=270,trim= 2cm 4cm 2cm 2cm]{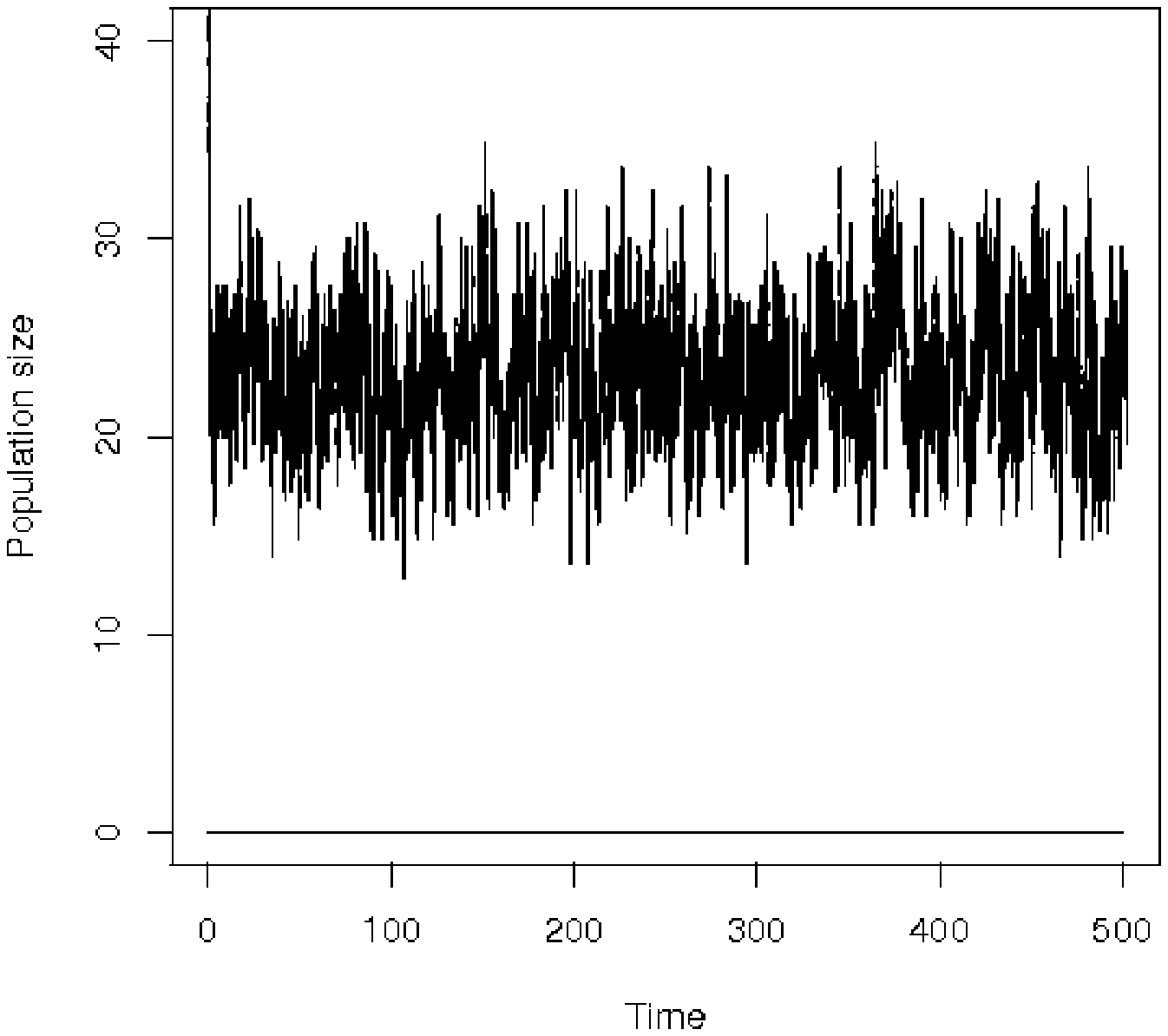} \vspace{0.2cm}&
\hspace{1.5cm}
\includegraphics[width=0.25\textwidth,height=0.15\textheight,angle=270,trim= 2cm 4cm 2cm 2cm]{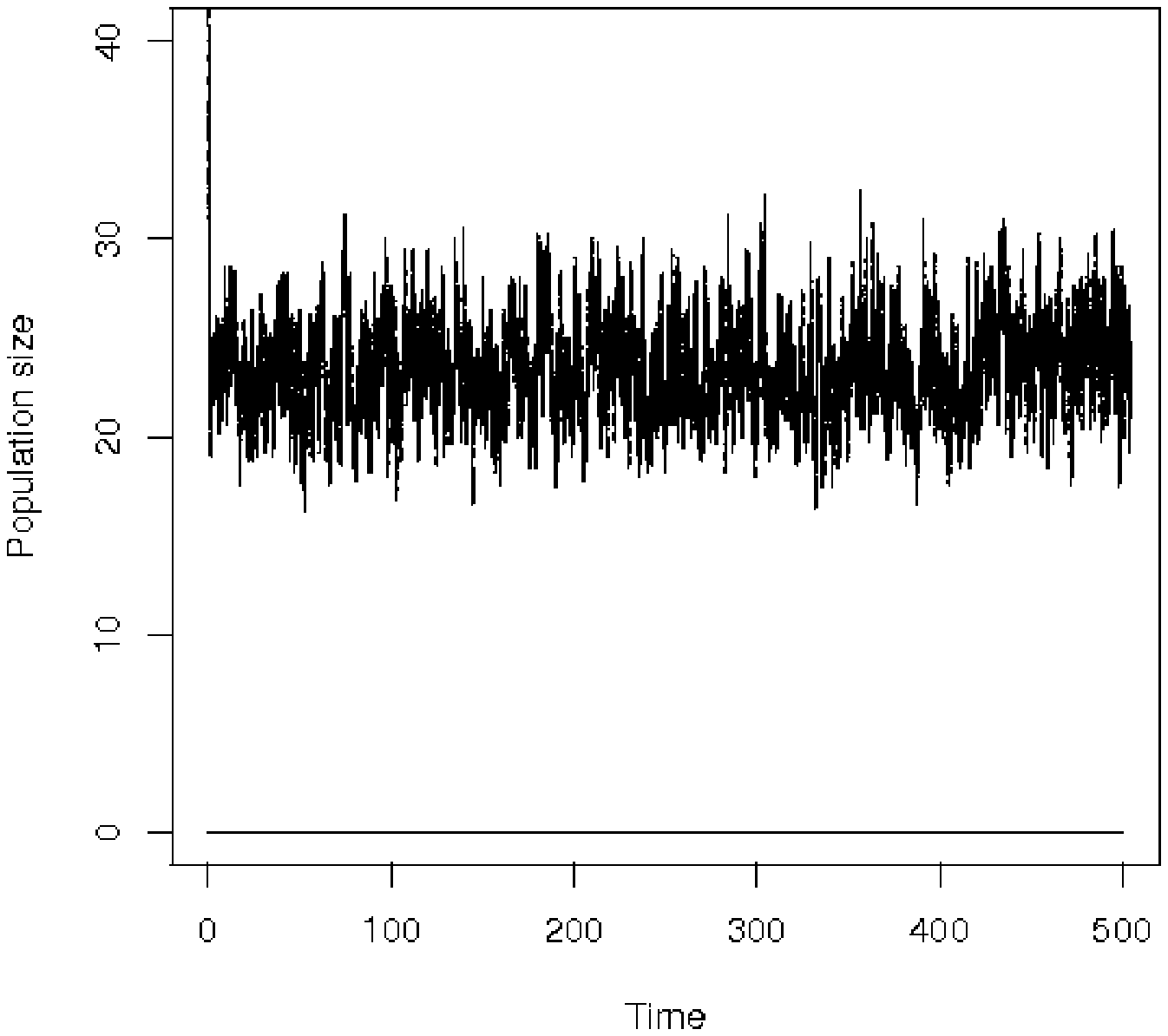}\vspace{0.2cm}\\
(a) & (b) & (c) \\
\end{tabular}
\end{center}
\vspace{-0.2cm}\caption{{\small \textit{Simulated size of a logistic age-structured population. The model is given by (\ref{EDS}) rescaled with $n=1$, $5$ and $10$ (from left to right).}}}\label{figuretaillepop}
\end{figure}

Even if simulations (shown in Figure \ref{figuretaillepop}) may lead us to the false impression that the long time behavior of the stochastic processes is the same as for their deterministic limit, the situation is much more complicated. We show that the stochastic process goes extinct almost surely, but after having spent an exponentially long time in the neighborhood of the stationary solution $\widehat{\xi}(da)=\widehat{n}(a)da$ (\ref{nchap}) of its deterministic approximation. Here, deterministic and stochastic models yield different conclusions, but the results for the limiting PDE make it possible to understand the behavior of the stochastic processes. With small system size the population size demonstrates wild fluctuations, and extinction is likely to occur quickly. As system size increases, fluctuations in population size are reduced, extinction takes more time, and the mean size of non-extinct populations converges toward the equilibrium solution of the deterministic model approximation.\\

\par First of all, on compact time intervals, the stochastic process belongs with large probability to a tube centered on the deterministic limit. This implies by Proposition \ref{proplgn} that:
\begin{equation}
\forall \varepsilon>0,\,\exists t_\varepsilon>0,\, \forall t>t_\varepsilon,\,\lim_{n\rightarrow +\infty}\mathbb{P}\left(\rho(Z^n_{t},\widehat{\xi})>\varepsilon\right)=0,\end{equation}where $\rho$ is for instance the Dudley metric which meterizes the topology of weak convergence on $\mathcal{M}_F(\widetilde{\X})$ (\eg \cite{rachev} p79). Heuristically, stochastic and deterministic processes have the same behavior on compact time intervals. In particular, since the deterministic process enters a given $\varepsilon$-neighborhood of its nontrivial stationary solution $\widehat{\xi}$ after a sufficiently large time that does not depend on $n$, the stochastic process enters an $2\varepsilon$-neighborhood of $\widehat{\xi}$ with a probability that tends to 1 when $n\rightarrow +\infty$. However, when $T\rightarrow +\infty$, the long time behavior differ.
\begin{prpstn}
For fixed $n$, one has almost sure extinction:
\begin{equation}
\mathbb{P}\left(\exists t\in \R_+,\, \langle Z^n_t,1\rangle =0\right)=1.
\end{equation}
\end{prpstn}Because of logistic competition, the population stays finite with a size that is controlled and cannot grow to infinity (see Proposition 5.6 \cite{trangdesdev}). The stochastic process finally leaves the neighborhood of $\widehat{\xi}$ to drive the population to extinction. We can thus ask how long the stochastic process stays in the neighborhood of $\widehat{\xi}$.
\par Let $0<\gamma<\langle \widehat{\xi},1\rangle$ and assume that we start our microscopic process in the neighborhood of the equilibrium $\widehat{\xi}$: $\rho(Z^n_0,\widehat{\xi})<\gamma$ a.s. An estimate of the time $\mathcal{T}^n=\inf\{t\in \R_+,\, \rho(Z^n_t,\widehat{\xi})\geq \gamma\}$ that the stochastic process spends in the neighborhood of $\widehat{\xi}$ can be established by large deviation results \textit{à la} Freidlin and Ventzell \cite{freidlinventzell, dembozeitouni} and generalized to our measure-valued setting (see Section 5.2 \cite{trangdesdev}). The proof (see \cite{trangdesdev}) is based on the fine properties of the convergence of $Z^n$ to its deterministic limit.
\begin{prpstn}$\exists \bar{V}>0,\,\underline{V}>0,\, \forall \delta>0,\,$
\begin{align}
\lim_{n\rightarrow +\infty}  \mathbb{P}(e^{n(\underline{V}-\delta)}<\mathcal{T}^n< e^{n(\bar{V}+\delta)})=1\label{estimeetn}
\end{align}
\end{prpstn}
When $n\rightarrow +\infty$, $\mathcal{T}^n$ tends to infinity in probability and extinction eventually disappears in the limit: unlike the microscopic process that gets extinct almost surely, its deterministic limit admits a nontrivial equilibrium. As a corollary, for every sequence $(t_n)_{n\in \N^*}$ such that $\lim_{n\rightarrow +\infty}t_n e^{-n(\underline{V}-\delta)}=0$ and $\lim_{n\rightarrow +\infty}t_n=+\infty$,
\begin{equation}
\lim_{n\rightarrow +\infty}\mathbb{P}\big(\rho(Z^n_{t_n},\widehat{\xi})>\varepsilon\big)=0.
\end{equation}
 Equation (\ref{estimeetn}) is interesting in itself to understand the persistence time of a stochastic finite population ; it is useful also to separate timescales in models of adaptive dynamics, as we shall see in the next sub-section.

\subsection{Application to the adaptive dynamics theory of phenotype evolution}\label{sectionAD}

In the theory of phenotype evolution, mutations are responsible for generating trait diversity in the population while natural selection acting through interaction and competition determines the traits that become fixed. Important questions in evolutionary biology require to include age structure in population models, as Charlesworth \cite{charlesworth}, Medawar \cite{medawar} and Stearns \cite{stearns} emphasized. The selective pressure and fixation probability of a trait may be functions of the age at which the trait is expressed, as well as of the age structure of the population. Reciprocally, the traits that are fixed can modify the age structure of the population.

\par The so-called "adaptive dynamics" framework considers large population limits and \textit{rare mutations}. Adaptive dynamics theory describes the evolution of traits in the population when it is possible to separate the timescales of ecology (birth and death events) and evolution (generation of new traits).
For populations with only trait structure, adaptive dynamics theory has been introduced and developed by Hofbauer and Sigmund \cite{hofbauersigmund}, Marrow \textit{et al.} \cite{marrow}, Dieckmann and Law \cite{dieckmannlaw}, Metz \textit{et al.} \cite{metzgeritzmeszenajacobsheerwaarden,durinxmetzmeszena} and Champagnat \etal  \cite{champagnat3, champagnatferrieremeleard, champagnatferrieremeleard2}.
\par These studies have been generalized to populations with age and trait structures by Méléard and Tran \cite{meleardtran}. Therein, equations of adaptive dynamics are derived from microscopic models by following Champagnat \cite{champagnat3} and several examples are studied and illustrated with numerical simulations. When mutations are sufficiently rare, natural selection wipes out the weakest competitors so that the resident population at the next mutation is monomorphic and in a stationary state. The estimates of extinction time provided in Section \ref{sectionextinction} tell us that mutations should not be \textit{too} rare if we wish to neglect the events of extinctions before the occurrence of a new mutant. By neglecting the transition periods, it is possible to describe evolution as a jump process by considering only the successive monomorphic equilibria (trait value and corresponding stationary age structure). If there exists a unique nontrivial stable stationary solution $\widehat{\xi}_x(da)=\widehat{n}(x,a)da$, it is sufficient to consider the trait valued process that jumps from one equilibrium trait to another. This process, called \textit{trait substitution sequence process for age-structured populations} (TSSPASP), introduced and analyzed in \cite{meleardtran}, can be described as follows.
 \begin{itemize}\item In a monomorphic population of trait $x$ at equilibrium, mutants are generated with rate $$p\,\int_{\R_+}b(x,a)\widehat{n}(x,a)da=p\widehat{n}(x,0)=\frac{p\widehat{N}_x}{\mathbb{E}(T_x)}$$where $\widehat{N}_x=\int_0^{+\infty}\widehat{n}(x,a)da$ is the size of the population at equilibrium and $\mathbb{E}(T_x)$ is the expected individual lifespan:$$\mathbb{E}(T_x)=\int_{0}^{+\infty} e^{-\int_0^a d(x,\alpha,\widehat{\xi}_xU(x,\alpha))d\alpha}da.$$
 \item The probability that the resulting mutant population with trait $x+h$ replaces the resident population with trait $x$ cannot always be calculated explicitly when we deal with populations that are structured by age. However, it can be computed numerically (\eg \cite{meleardtran}).
\end{itemize}
\par The theory of random point processes is necessary to model these events that occur at discrete random times. Evolutionary models have also been developed from a deterministic point of view, but with different renormalizations that do not describe the same \textit{rare} mutations (Diekmann \etal \cite{diekmannjabinmischlerperthame}, Carrillo \etal \cite{carrillocuadradoperthame}).

\par If mutation steps are additionally assumed small, the TSSPASP can be approximated by the solution of an ordinary differential equation that generalizes the "\textit{canonical equation of adaptive dynamics}" (see Dieckmann and Law \cite{dieckmannlaw}, Champagnat \cite{champagnat}).
\par Recently, models of populations with age structure have been considered in adaptive dynamics theory dealing with functional traits (\eg  \cite{dieckmannheinoparvinen,ernandedieckmanheino,parvinendieckmannheino}). The traits $x(a)$ are functions of age (growth curve, flowering intensity, birth or death curves...) In Dieckmann \etal \cite{dieckmannheinoparvinen} or Parvinen \etal \cite{parvinendieckmannheino}, the birth and death rates are functions of the functional trait, but not of age (for example, the averaged form $b(x)=B(\int_0^{+\infty}x(a)da)$ with $B\in \Co_b(\R,\R_+)$ was assumed).

\section{Age and trait structures: model examples and insights from simulations}\label{sectionsimu}

We now present several examples of interesting biological questions involving populations with age and trait structures. Mechanisms are described at the individual level, and we are interested in the resulting macroscopic dynamics. Not surprisingly, the more realistic and complex the model is, the more difficult it is to carry out an analytical study of the model behavior. Simulations offer an alternative and complementary approach.

\subsection{Simulation algorithm}\label{sectionsimulationsage}

\par Following the algorithms proposed by Fournier and Méléard \cite{fourniermeleard}, it is possible to simulate exactly the law of the process $Z$ (\ref{EDS}), without approximation scheme or grid. The main difficulty comes from the fact that the rates depend on age that changes continuously in time. To handle this, we use acceptance-reject procedures.
\par Let $Z_0\in \mathcal{M}_P(\widetilde{\X})$ be an initial condition. We simulate by recursion a succession of birth and death events that modify the size of the population. Let us set $T_0=0$ for the event number 0. Assume that we have already simulated $k$ events ($k\in \N$), and that the last of these events had occurred at time $T_k$. The size of the population $N_{T_k}=\langle Z_{T_k},1\rangle $ at time $T_k$ is finite and the global jump rate is upper bounded by $\bar{b} N_{T_k}+\bar{d}(1+N_{T_k})N_{T_k}$, which is finite. To obtain $T_{k+1}$ for the $k+1^{\mbox{th}}$ event, we simulate candidate events from time $T_k$. The latter are given by a sequence $(\tau_{k,\ell})_{\ell\in \N}$ of possible event times following a Poisson point process with intensity $\bar{b} N_{T_k}+\bar{d}(1+N_{T_k})N_{T_k}$. The first of these times which is accepted by the procedure defines $T_{k+1}$.
\begin{enumerate}
\item[0.] We set $\tau_{k,0}:=T_k$, $N_{\tau_k}:=N_{T_k}$ and $\ell:=0$.
\item[1.] We simulate independent exponential variables $\varepsilon_{k,\ell}$ with parameter 1 and we define $\tau_{k,\ell+1} = \tau_{k,\ell} +
\varepsilon_{k,\ell}/[\bar{b}N_{T_{k}}+\bar{d}(1+N_{T_{k}})N_{T_k}]$.
\item[2.] On the interval $[\tau_{k,\ell},\, \tau_{k,\ell+1}[$, only aging takes place.
 For $i\in \lbrac 1,N_{T_k}\rbrac$, the age of individual $i$ becomes $A_i(\tau_{k,\ell})+\tau_{k,\ell+1}-\tau_{k,\ell}$ at time $\tau_{k,\ell+1}$.
\item[3.] We simulate an integer valued r.v. $I_{k,\ell}$ uniformly distributed on
$\lbrac 1,N_{T_{k}}\rbrac$, and we define the following quantities in $[0,1]$:{\small
\begin{align}
 & \widetilde{m}_1(\tau_{k,\ell+1},Z_{\tau_{(k,\ell+1)_-}},I_{k,\ell})  =
\frac{\big(1-p\big) b(X_{I_{k,\ell}}(\tau_{(k,\ell+1)-}),A_{I_{k,\ell}}(\tau_{(k,\ell+1)_-}))}{\bar{b}+\bar{d}(1+N_{T_{k}})}\nonumber\\
 & \widetilde{m}_2(\tau_{k,\ell+1},Z_{\tau_{(k,\ell+1)_-}},I_{k,\ell})
 =    \widetilde{m}_1(\tau_{k,\ell+1},Z_{\tau_{(k,\ell+1)_-}},I_{k,\ell})  +  \frac{p\,b(X_{I_{k,\ell}}(\tau_{(k,\ell+1)-}),A_{I_{k,\ell}}(\tau_{(k,\ell+1)_-}))}{\bar{b}+\bar{d}(1+N_{T_{k}})}\nonumber\\
  & \widetilde{m}_3(\tau_{k,\ell+1},Z_{\tau_{(k,\ell+1)_-}},I_{k,\ell})  =  \widetilde{m}_2(\tau_{k,\ell+1},Z_{\tau_{(k,\ell+1)_-}},I_{k,\ell}) \nonumber\\
  & +  \frac{d\left(X_{I_{k,\ell}}(\tau_{(k,\ell+1)-}),A_{I_{k,\ell}}(\tau_{(k,\ell+1)_-}),
Z_{\tau_{(k,\ell+1)_-}}U(X_{I_{k,\ell}}(\tau_{(k,\ell+1)-}),A_{I_{k,\ell}}(\tau_{(k,\ell+1)_-}))\right)}{
\bar{b}+\bar{d}(1+N_{T_{k}})}\nonumber
\end{align}}
\item[4.] We simulate a r.v. $\Theta_{k,\ell}$ with uniform law on $[0,1]$.
\begin{enumerate}
\item If $0\leq \Theta_{k,\ell} <\widetilde{m}_1(\tau_{k,\ell+1},Z_{\tau_{(k,\ell+1)_-}}, I_{k,\ell})$, then the $k+1^{\mbox{th}}$ event occurs: Individual $I_{k,\ell}$ gives birth to a clone with age $0\in \R_+$ and one defines $T_{k+1}=\tau_{k,\ell+1}$.
\item If $\widetilde{m}_1(\tau_{k,\ell+1},Z_{\tau_{(k,\ell+1)-}}, I_{k,\ell})\leq \Theta_{k,\ell} <\widetilde{m}_2(\tau_{k,\ell+1},Z_{\tau_{(k,\ell+1)_-}}, I_{k,\ell})$, then the $k+1^{\mbox{th}}$ event occurs: Individual $I_{k,\ell}$ gives birth to a mutant of age 0 and trait $x+h$ where $h$ follows the distribution probability of density $k(X_{I_{k,\ell}}(\tau_{(k,\ell+1)-}),A_{I_{k,\ell}}(\tau_{(k,\ell+1)-}),h')$. One defines $T_{k+1}=\tau_{k,\ell+1}$.
    \item If
$\widetilde{m}_2(\tau_{k,\ell+1},Z_{\tau_{(k,\ell+1)-}}, I_{k,\ell})\leq \Theta_{k,\ell} <\widetilde{m}_3(\tau_{k,\ell+1},Z_{\tau_{(k,\ell+1)_-}}, I_{k,\ell})$, then the $k+1^{\mbox{th}}$ event occurs: Individual $I_{k,\ell}$ dies and one defines $T_{k+1}=\tau_{k,\ell+1}$.
\item If $\widetilde{m}_3(\tau_{k,\ell+1},Z_{\tau_{(k,\ell+1)-}}, I_{k,\ell})\leq \Theta_{k,\ell}$ then nothing happens. We reiterate the algorithm from $1.$ with $\ell+1$ in place of $\ell$ until we obtain the $k+1^{\mbox{th}}$ event.
\end{enumerate}
\end{enumerate}

\begin{rmrk}
To each individual are associated three clocks corresponding respectively to birth (with or without mutation) and death events. The event that actually occurs corresponds to the minimum of these durations, over the set of all individuals.
\end{rmrk}

\subsection{Examples}\label{sectionsimulationsexemple}

\par Simulations of particular cases are now presented. For each example, simulations have been run several times ; only one particular simulation is presented for each example as an illustration. The programs have been written with the \textbf{R} freeware (\footnote{http://www.r-project.org/}) and \textsc{Matlab}.

\subsubsection{Example 1: Evolution of offspring size in age structured populations with size-dependent competition}\label{sectionexemple1}

\par In this example, we are interested in a population with size-dependent competition, and ask how the evolution of body size changes when we take the growth of individuals into account.
\par Individuals are characterized by their body size at birth $x_0\in [0,4]$ which is the heritable trait subject to mutation, and by their physical age $a\in [0,2]$, with aging velocity $1$. Body size is an increasing function of age:
\begin{equation}x=x_0+g\, a,\label{tailleexemple1}\end{equation}where $g$ is the growth rate, which is assumed constant and identical for all individuals.
\par An individual of trait $x_0\in [0,4]$ gives birth at the age independent rate:
\begin{equation}\label{simulationtauxfecondite}
b(x_0)=4-x_0.\end{equation}This is a simple expression of the standard trade-off between fecundity and offspring size (\eg Stearns \cite{stearns}). With probability $p=0.03$ a mutation occurs and affects $x_0$. The new trait is $x'_0  =   \min ( \max (0, x_0+y_{x_0}) , 4)$, where $y_{x_0}$ is a Gaussian r.v. with expectation 0 and variance $0.01$. With probability $1-p$, the offspring inherits its parent's trait, $x_0$.

\par The death rate of an individual with trait $x_0\in [0,4]$ and age $a\in [0,2]$ living in a population $Z\in \mathcal{M}_{P}([0,4]\times [0,2])$ is given by:
\begin{align}
d(x_0,a,Z) = & \int_{\widetilde{\X}}U (x_0+g\, a - x'_0-g\, \alpha)Z(dx'_0,d\alpha),\label{simulationtauxmort}
\end{align}where\begin{align}U(x-y)= & \frac{2}{300}\left(1-\frac{1}{1+1.2\,\exp\left(-4(x-y)\right)}\right)\in \left[0,\frac{2}{300}\right]\label{noyaukisdi}\end{align}
is the asymmetric competition kernel introduced by Kisdi \cite{kisdi} which gives a competitive advantages to larger individuals. The important point is that the size of competitors, hence the intensity of competition between them, are not constant when $g\not=0$. Even in a monomorphic population, size varies across individuals and competition experienced by each individual varies with its size, which is age-dependent.

\begin{figure}[!ht]
\begin{center}
\begin{tabular}[!ht]{ccc}
\includegraphics[width=0.25\textwidth,height=0.15\textheight,angle=270,trim= 2cm 4cm 2cm 2cm]{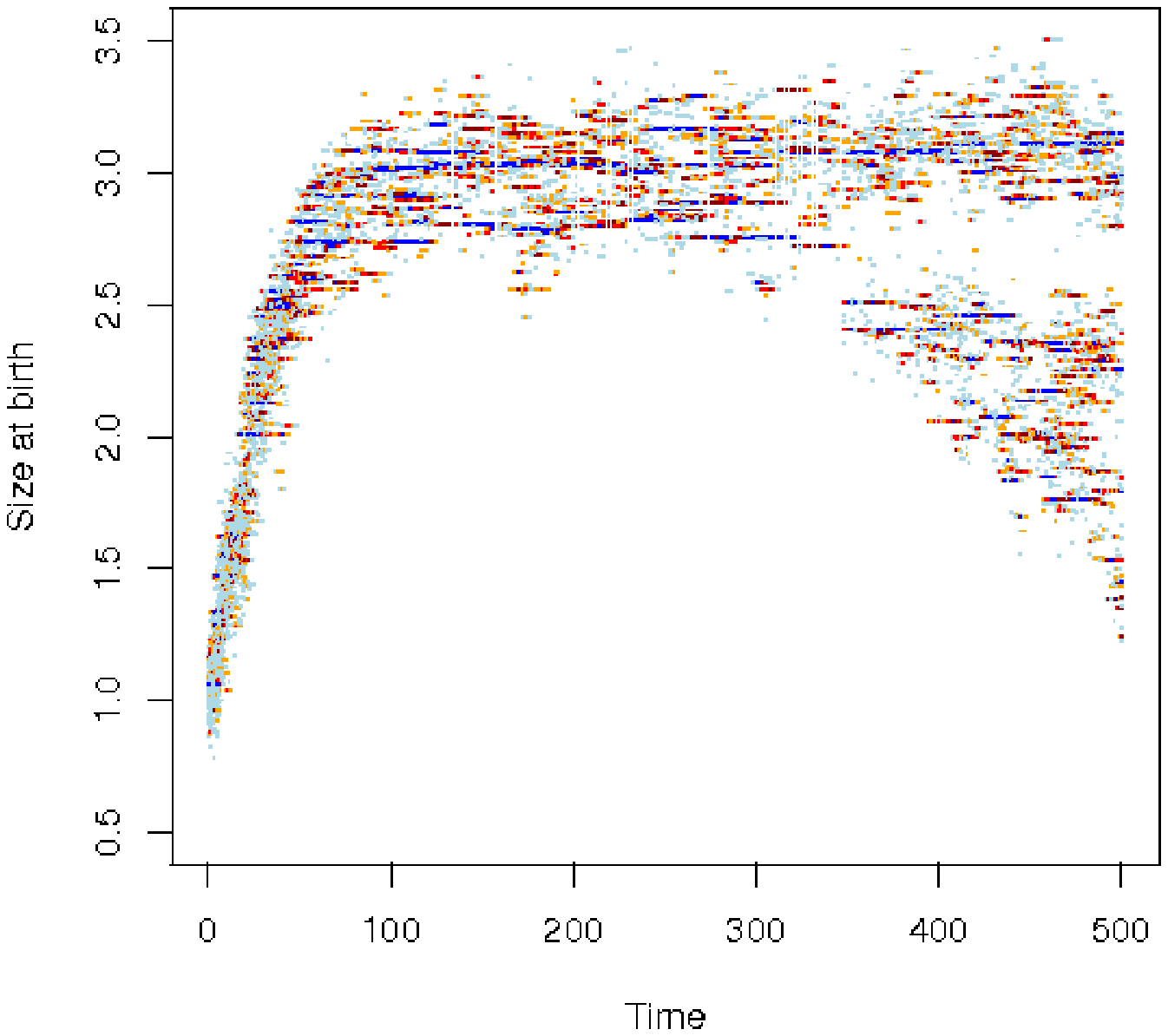}
\vspace{0.25cm}
&
\hspace{2cm}
\includegraphics[width=0.25\textwidth,height=0.15\textheight,angle=270,trim= 2cm 4cm 2cm 2cm]{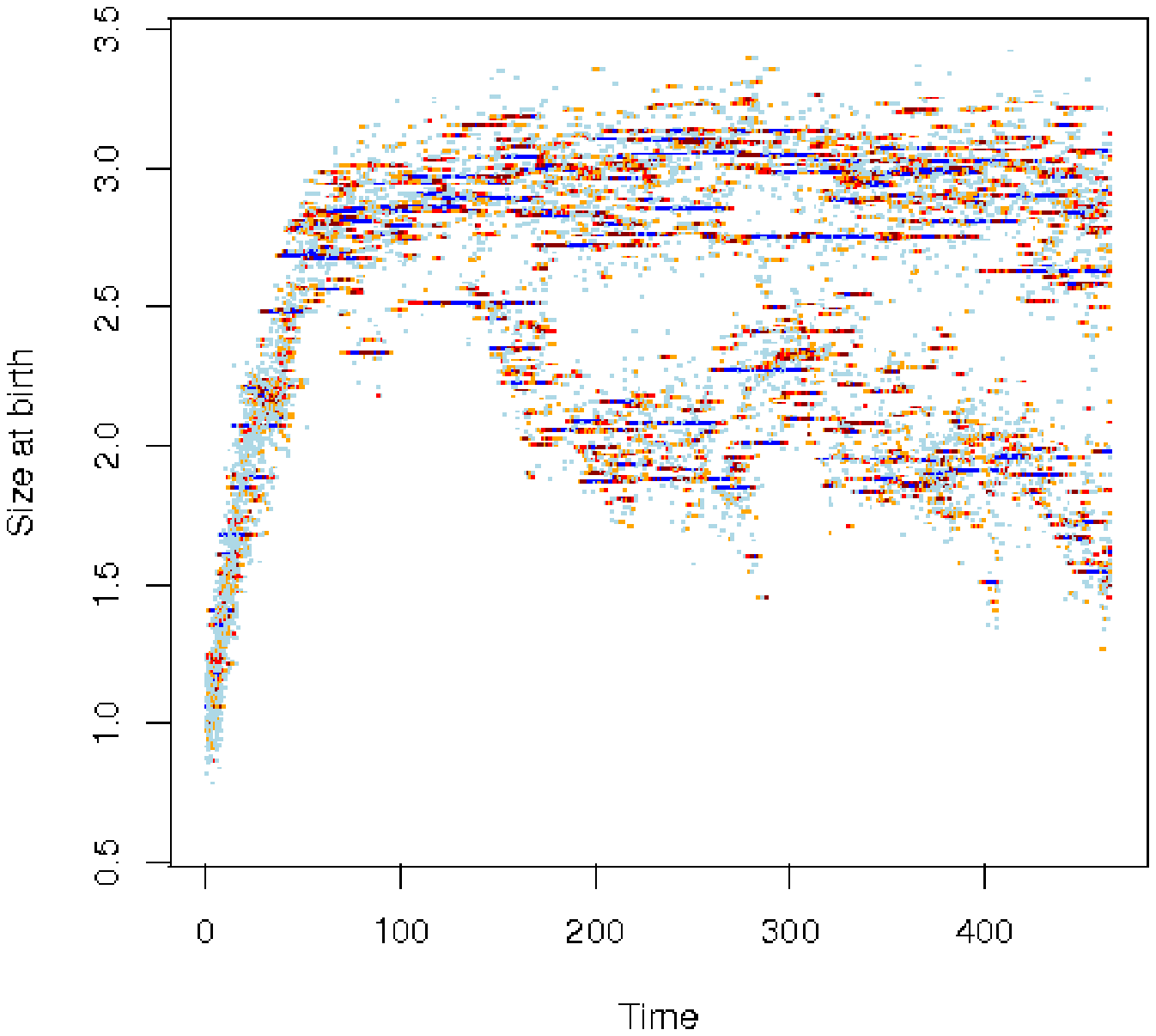}
\vspace{0.25cm}
&
\hspace{2cm}
\includegraphics[width=0.25\textwidth,height=0.15\textheight,angle=270,trim= 2cm 4cm 2cm 2cm]{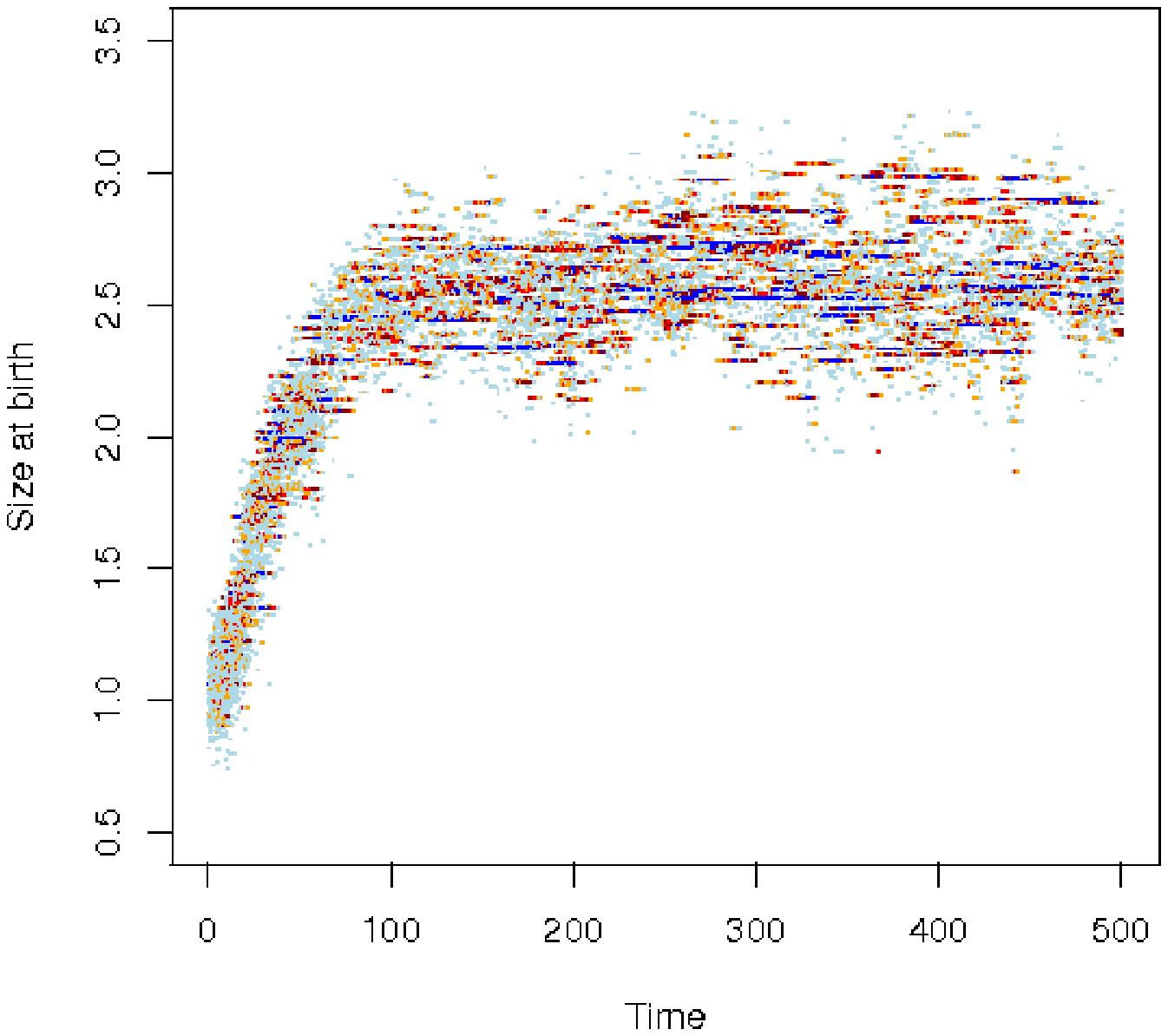}
\vspace{0.25cm}
\\
(a) &
(b) &
(c)
\end{tabular}
\end{center}
\vspace{-0.2cm}\caption{{\small \textit{Trait dynamics in an age structured population with size-dependent competition. The evolving trait is offspring size. See Example 1 in main text for details. (a) $g=0$, (b) $g=0.3$, (c) $g=1$.}}}\label{figsim1}
\vspace{-0.5cm}
\end{figure}

\par In the simulations of Figure \ref{figsim1}, we chose $g=0$, $g=0.3$ and $g=1$. The population at time $t=0$ is monomorphic with trait $x_0=1,06$, and contains $N_0=900$ individuals. The initial age distribution is uniform on $[0,2]$, and the corresponding sizes are computed using (\ref{tailleexemple1}). For each $t$ the support of the measure $\sum_{i=1}^{N_t} \delta_{(x_{0})_i}$ is represented.
\par For $g=0$ (constant size during life), we observe a branching phenomenon around $t= 300$: the distribution of trait $x_0$ splits into two "branches". The initial branch first stabilizes around some equilibrium ($x_0\simeq 2.7$), and then, a second branch appears below ($x_0\simeq 2$). For $g=0.3$, the population branches into two subpopulations ($x_0\simeq 2.7$ and $x_0\simeq 2$) sooner, around $t=100$. For $g=1$, the branching phenomenon is not observed anymore. The occurrence of branching in Fig.\ref{figsim1} (a) shows that the coexistence of subpopulations can be promoted by the competition-fecundity trade-off: small individuals may compensate for their competitive inferiority by reproducing more often. In (b), non-zero body growth results in a wider size-structure in the population, so that the conditions of branching are met earlier. Here, growth acts as a mixing factor, which lessens the differences between traits. In (c), individual growth is so fast that the differentiation of a subpopulation of smaller individuals is always prevented by strong competition with very large individuals.

\subsubsection{Example 2: Coevolution of offspring size, growth rate and age at maturity under size and stage-dependent competition}\label{sectionsimulationsoscillations}

As in Example 1, we consider an age- and trait-structured population where the traits are the size at birth $x_0\in [0,4]$, the growth rate $g\in [0,2]$, and the age at maturity $a_M\in [0,2]$. The scalar age $a\in [0,2]$ (with aging velocity 1) still corresponds to the physical age.
\par The age at maturity $a_M$ separates the life history of individuals into two periods: a growth period and a reproduction period. Calsina and Cuadrado \cite{calsinacuadrado1} also studied a population that is structured into two age classes, juveniles and adults. Assuming that the length of the juvenile period is exponentially distributed, they described population dynamics by a system of two ODEs. They then modelled the evolution of the age at maturity by considering whether a small mutant population can invade the resident population, leading to questions of stability of the stationary solutions. Ernande \textit{et al.} \cite{ernandedieckmanheino} considered a similar problem: their population was divided into age classes, and the age at maturity was plastic and determined by a heritable reaction norm that was subject to mutation.

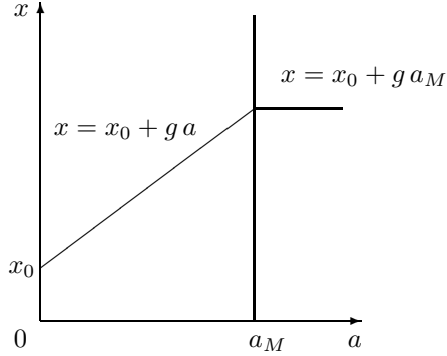
\begin{figure}[!ht]
\begin{center}
    \begin{picture}(150,130)(-20,-10)
      \put(0,0){\vector(1,0){120}} \put(0,0){\vector(0,1){120}}
      \put(-10,-10){0}
      \put(78, -10){$a_M$}
       \put(115,-10){$a$}
        \put(-10, 115){$x$}
        \put(-12, 18){$x_0$}
        \put(5, 70){$x=x_0+g\,a$}
        \put(90, 90){$x=x_0+g\,a_M$}
      \put(80,0){\line(0,1){115}}
      \put(0,20){\line(4,3){80}}
      \put(80,80){\line(1,0){33}}
    \end{picture}
    \end{center}
\caption{{\small \textit{Size $x$ as a function of age $a$. The age at maturity is denoted by $a_M$. Before maturity, juveniles born with size $x_0$ grow at rate $g$. After maturity, the size of adults is constant.}}}\label{fig3}
\end{figure}
In our model, the size of an individual with trait $(x_0,g,a_M)$ and age $a$ is given by:
\begin{equation}
x(a,x_0, g,a_M)=x_0+g\, (a\wedge a_M).\label{tailleexemple2}\end{equation}
Before maturity, individuals (called juveniles) do not reproduce and invest all their resources into growth. Once the age $a_M$ has been reached, individuals (called adults) have constant size and begin to reproduce (see Fig.\ref{fig3}).
\par The birth rate of an individual with traits $(x_0,g,a_M)\in [0,4]\times [0,2]\times [0,2]$ is 0 for juveniles and a decreasing function of $x_0$ for adults:
\begin{equation}
b(x_0,a_M,a)  =  \left(4-x_0\right)\mathbf{1}_{a\geq a_M}.
\end{equation}
The traits $x_0$, $g$ and $a_M$ are heritable and can be affected in offspring by mutation that occur with probability $p=0.03$. Mutated traits are chosen according to:
\begin{align*}
x'_0  =   \min ( \max (0, x_0+y_{x_0}) , 4)\quad g'  =   \min ( \max (0, g+y_g) , 2),\quad a'_M  =   \min ( \max (0, a_M+y_{a_M}) , 2),
\end{align*}where $y_{x_0}$, $y_g$ and $y_{a_M}$ are Gaussian independent r.v. with expectation equal to 0 and variance equal to $0.01$.

\par The death rate of an individual with traits $(x_0,g,a_M)\in [0,4]\times [0,2]\times [0,2]$ and age $a\in [0,2]$ in the population $Z\in \mathcal{M}_P([0,4]\times [0,2]\times [0,2]\times [0,2])$ is given by:
\begin{align*}
d(x_0 ,g,a_M,a,Z)
=  \left\{\begin{array}{l}d_0\times g\, +10^{-5} \int_{\widetilde{\X}}U(x_0+g\,a-x'_0-g'(a'\wedge a'_M))Z(dx'_0,dg',da'_M,da') \quad  \mbox{ if }a\leq a_M,\, \\
 \int_{\widetilde{\X}}U(x_0+g\,a_M-x'_0-g'(a'\wedge a'_M))Z(dx'_0,dg',da'_M,da') \quad \mbox{ if } a_M < a \leq 2,\, \\
+\infty \quad  \mbox{ if }a> 2,\end{array}\right.
\end{align*}where $U$ is the competition kernel introduced in (\ref{noyaukisdi}) and $d_0\geq 0$ measures the adversity of environmental conditions. During the juvenile period the main factor of death is related to growth: individuals with high growth rates $g$ need more ressources and incur a survival cost of foraging ; thus, the density-independent component $d_0 g$ of their death rate is higher if $d_0$ is larger. Then, during the adult period, there is no more growth and death is fully determined by logistic competition, which was soften for juvenile individual. This may apply to species that shift to different habitats or resources as individuals reach maturity.
\begin{figure}[!ht]
\begin{center}
\begin{tabular}[!ht]{cc}
(a)
\vspace{0.2cm} & (b)\vspace{0.2cm} \\
\includegraphics[width=0.18\textwidth,height=0.15\textheight,angle=270,trim= 2cm 4cm 2cm 2cm]{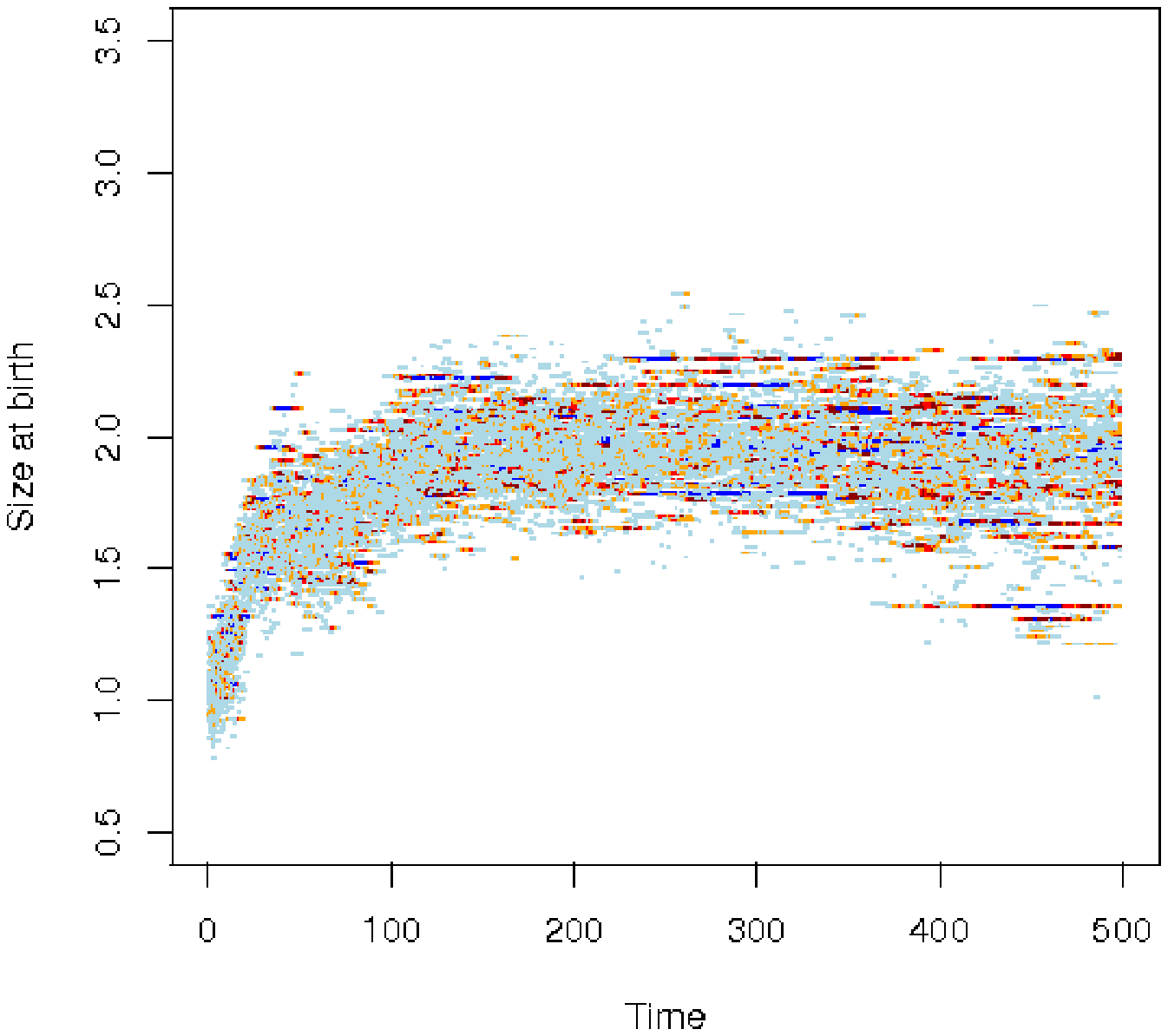}
\par\vspace{0.4cm}
 & \hspace{2cm}
\includegraphics[width=0.18\textwidth,height=0.15\textheight,angle=270,trim= 2cm 4cm 2cm 2cm]{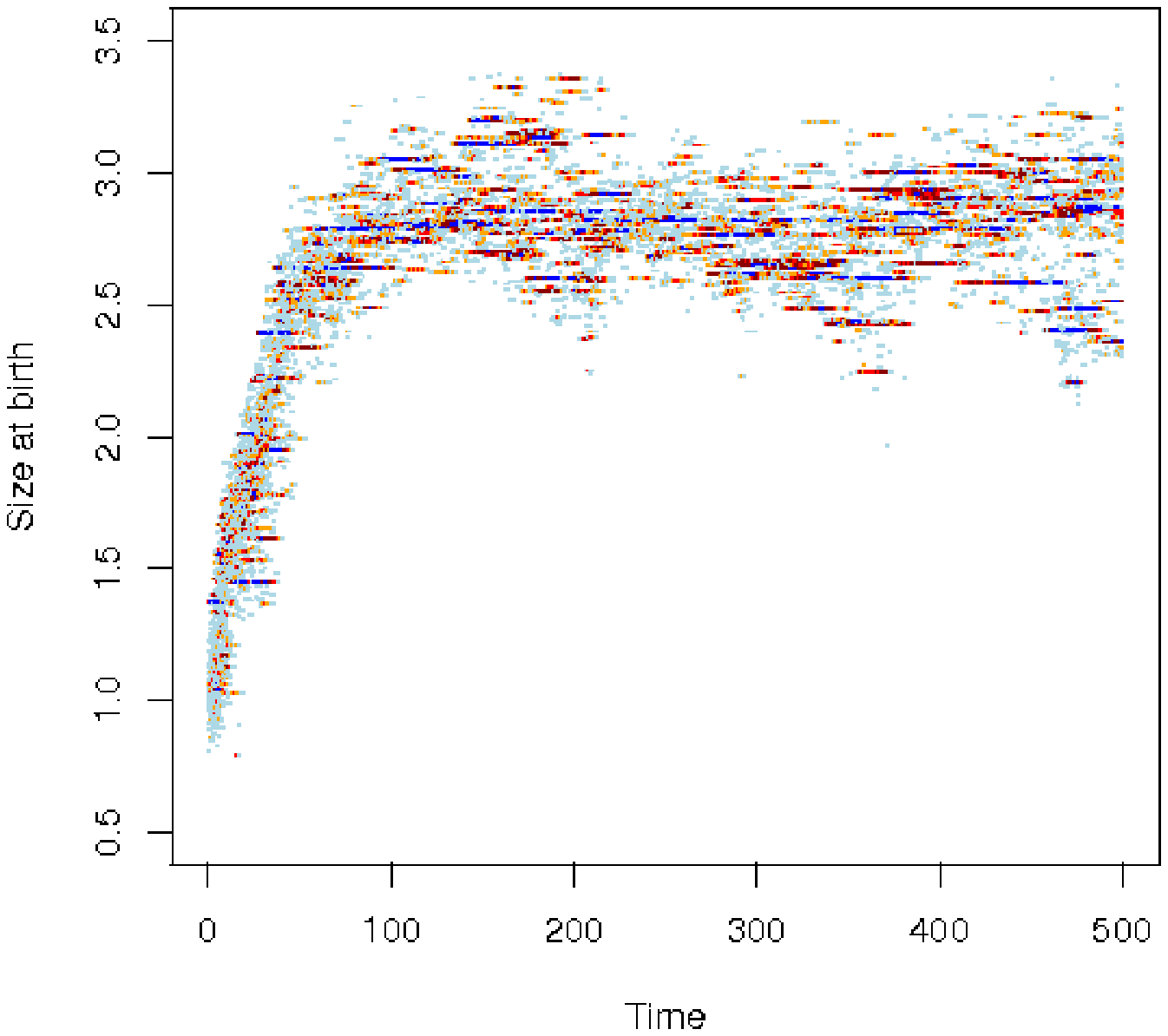}
\par\vspace{0.4cm}
  \\
\includegraphics[width=0.18\textwidth,height=0.15\textheight,angle=270,trim= 2cm 4cm 2cm 2cm]{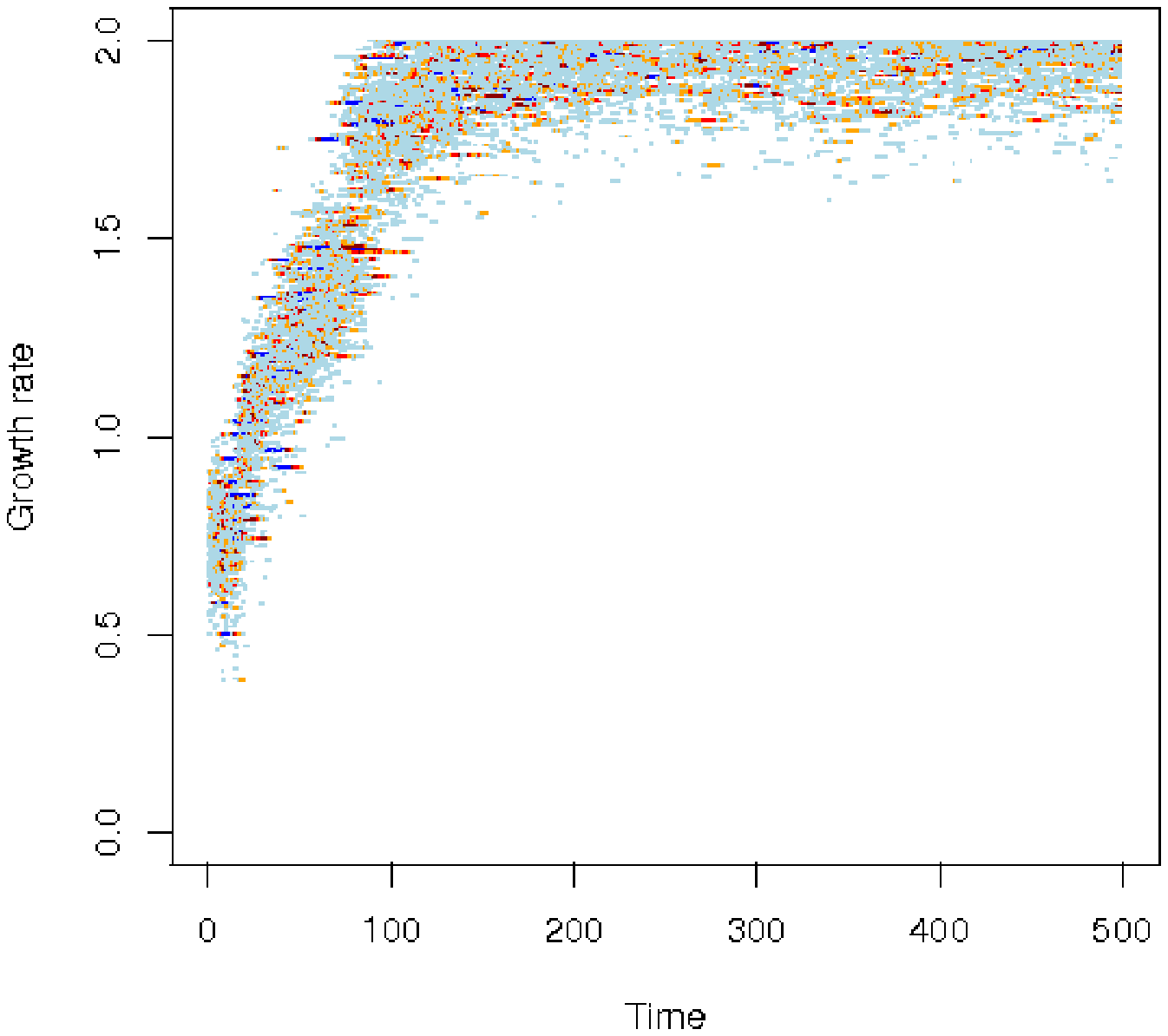}
\par\vspace{0.4cm}
 &
\hspace{2cm}\includegraphics[width=0.18\textwidth,height=0.15\textheight,angle=270,trim= 2cm 4cm 2cm 2cm]{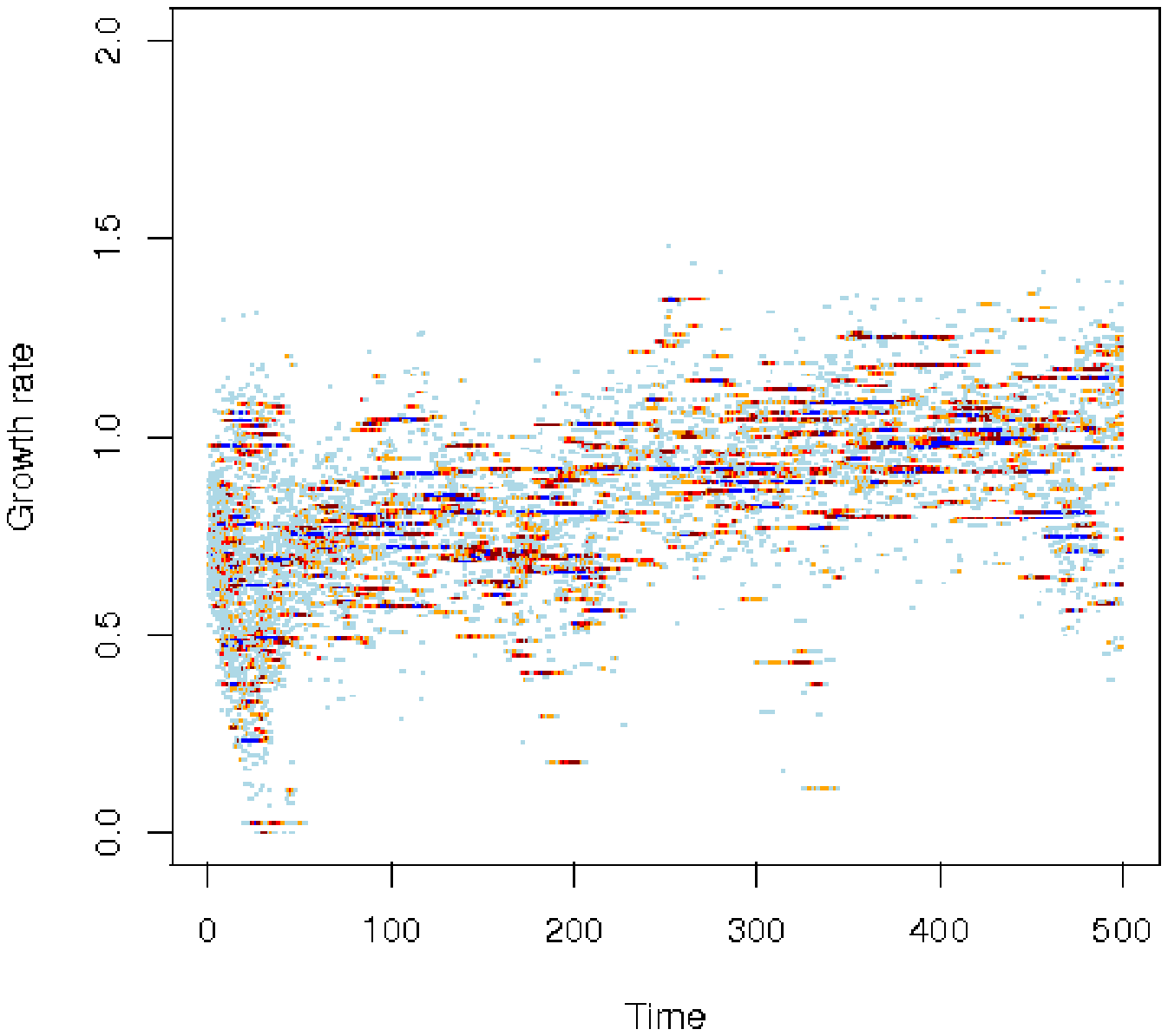}
\par\vspace{0.4cm}
 \\
\includegraphics[width=0.18\textwidth,height=0.15\textheight,angle=270,trim= 2cm 4cm 2cm 2cm]{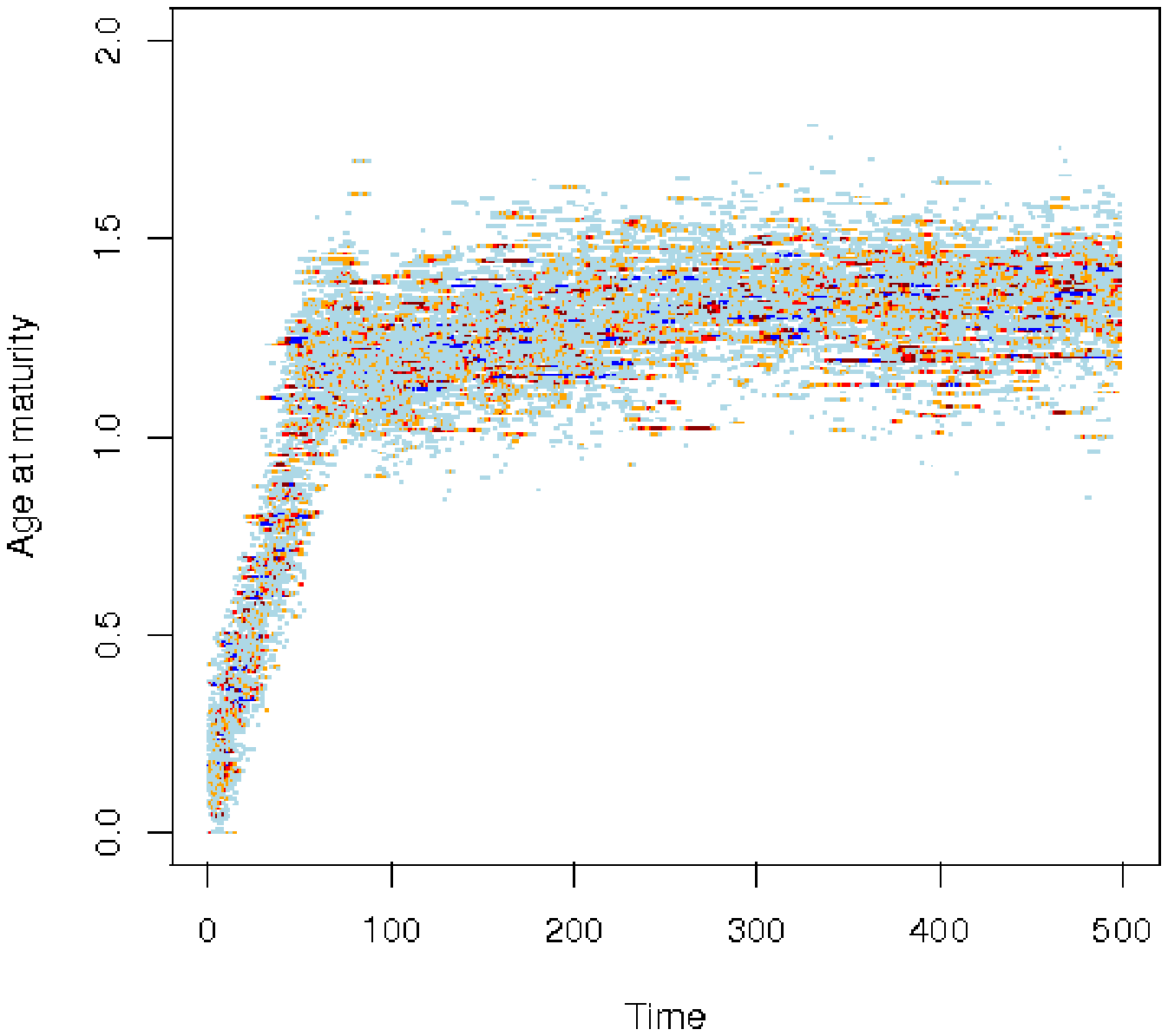}
\par\vspace{0.4cm}
 &
\hspace{2cm}\includegraphics[width=0.18\textwidth,height=0.15\textheight,angle=270,trim= 2cm 4cm 2cm 2cm]{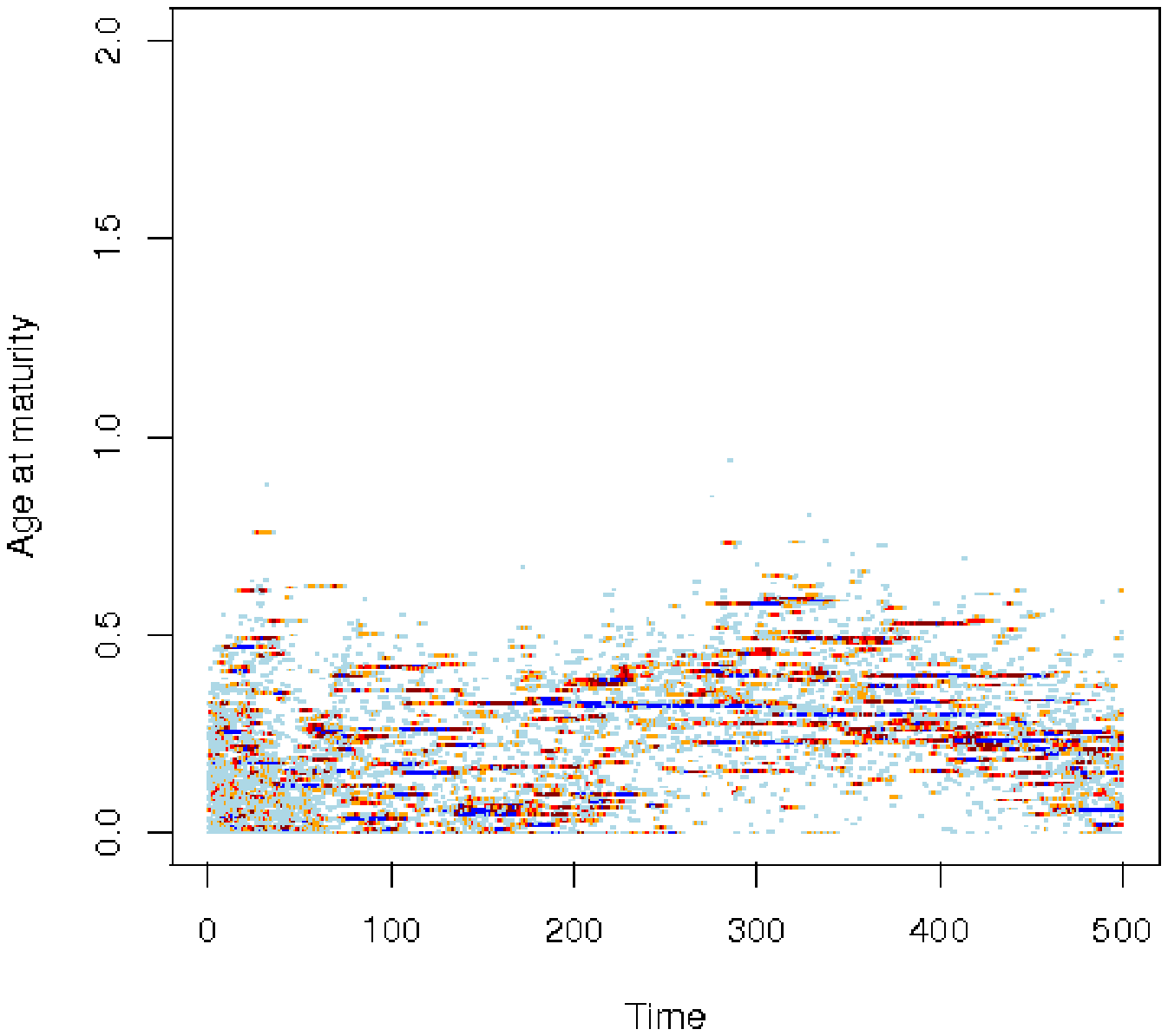}
\par\vspace{0.4cm}
 \\
\includegraphics[width=0.18\textwidth,height=0.15\textheight,angle=270,trim= 2cm 4cm 2cm 2cm]{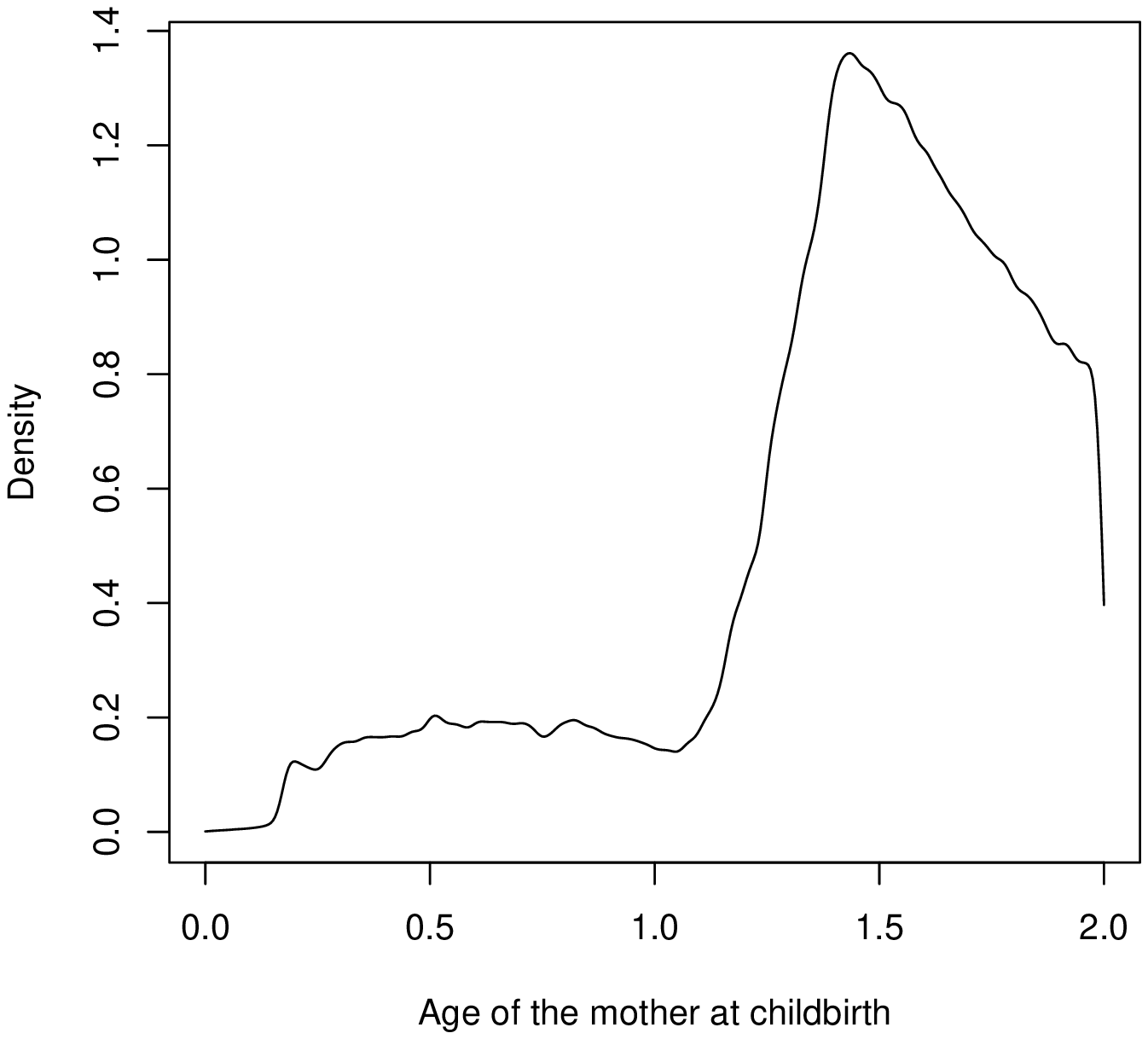}
\par\vspace{0.2cm}
 &
\hspace{2cm}\includegraphics[width=0.18\textwidth,height=0.15\textheight,angle=270,trim= 2cm 4cm 2cm 2cm]{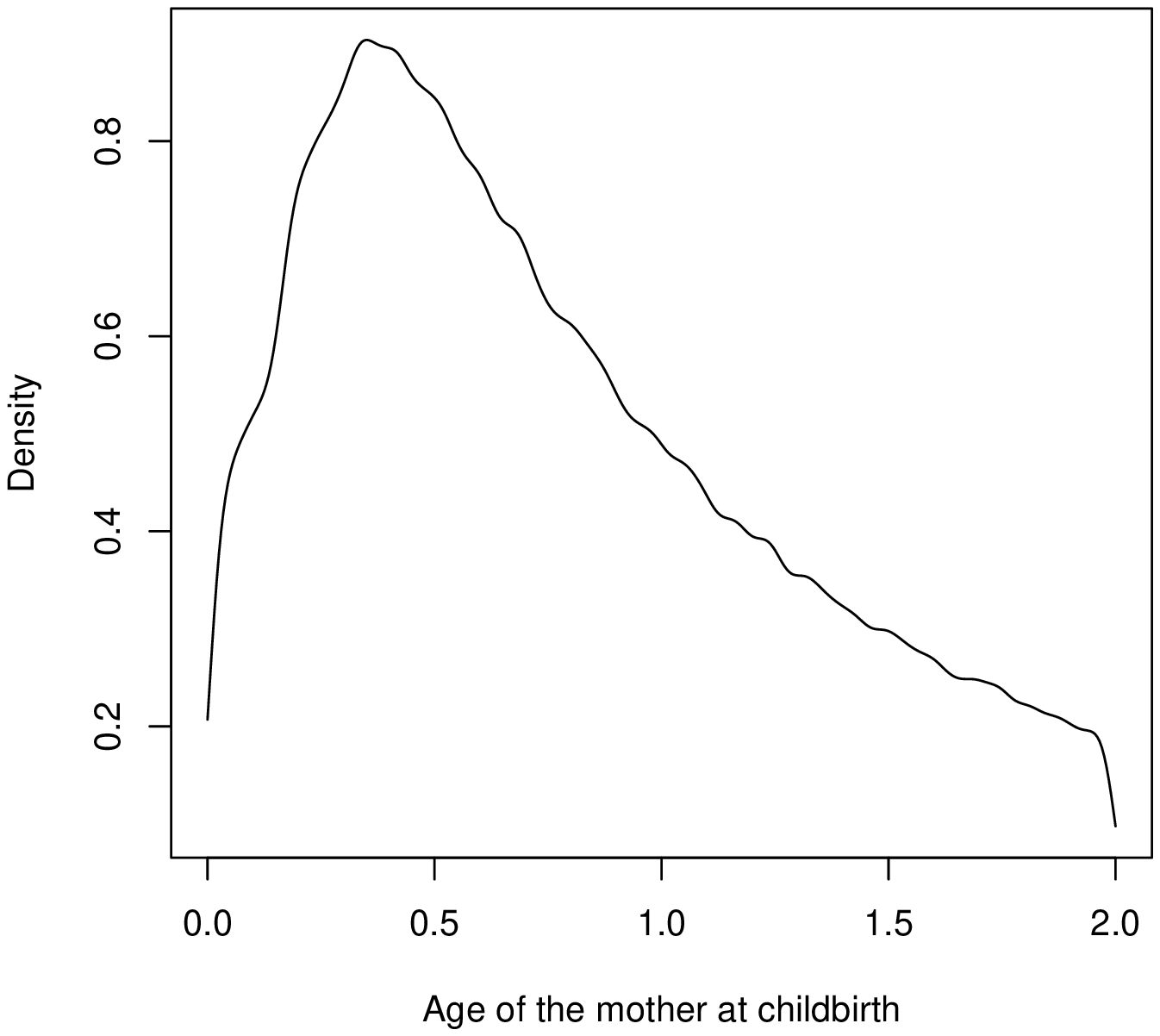}
\par\vspace{0.2cm}
\end{tabular}
\caption{{\small \textit{Time dynamics of three co-evolving traits: size at birth, growth rate, age at maturity (rows 1-3), and population structure (distribution of maternal ages at birth) in the "equilibrium state" (row 4). See Exemple 2 in main text for details. (a): $d_0=3.3\,10^{-3}$, (b): $d_0=0.5$.
}}}\label{figsim11}
\end{center}
\end{figure}
\par As in Example 1, we chose as initial state a monomorphic population ($x_0=1.06$, $g=0.74$ and $a_M=0.20$) with $N_0=900$ individuals whose age is uniformly drawn in $[0,2]$. Fig. \ref{figsim11} shows that changing $d_0$ can result in very different dynamics. When $d_0$ is small (Fig.\ref{figsim11}(a), $d_0=3.3\, 10^{-3}$), individuals with higher growth rates do not pay much of a cost ($\forall g\in [0,2],\, 0\leq d_0 g\leq 6.6\,10^{-3}$). The juvenile death rate is small (of the order of $10^{-2}$ for a population of 1000 individuals) compared to the adult death rate (of the order of 1). Since competition favors larger individuals (the term $U(x_0+ga_M-x)$ describing competition exerted by an individual with size $x$ on an adult with size $x_0+ga_M$ is smaller when $x_0+ga_M$ is large). Hence there is a strong selective advantage to grow fast ($g\simeq 2$) and long ($a_M \simeq 1.3$). The reproduction loss caused by delayed maturity is further compensated by the evolution of relatively small size at birth. Consequences for population dynamics are shown in Fig.\ref{figsim12}. The population stabilizes numerically around a state with high growth rate and late maturation.

\par
\begin{figure}[!ht]
\begin{center}
\begin{tabular}[!ht]{cc}
(a) & (b) \\
\hspace{1cm}\begin{minipage}[b]{.4\textwidth}\centering
\includegraphics[width=0.65\textwidth,height=0.22\textheight,angle=270,trim= 2cm 4cm 2cm 2cm]{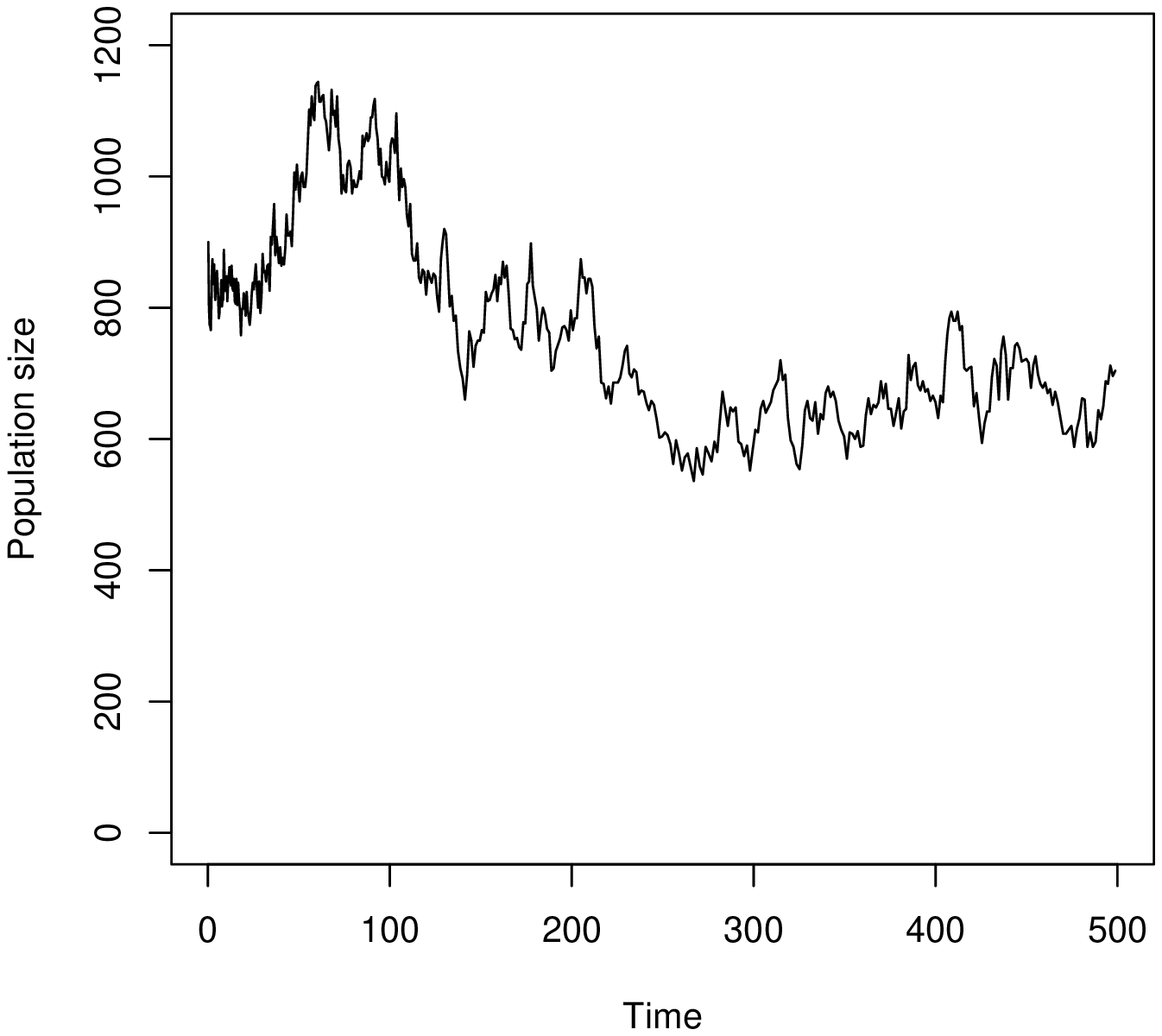}
\end{minipage} & \hspace{2cm}
\begin{minipage}[b]{.4\textwidth}\centering
\includegraphics[width=0.65\textwidth,height=0.22\textheight,angle=270,trim= 2cm 4cm 2cm 2cm]{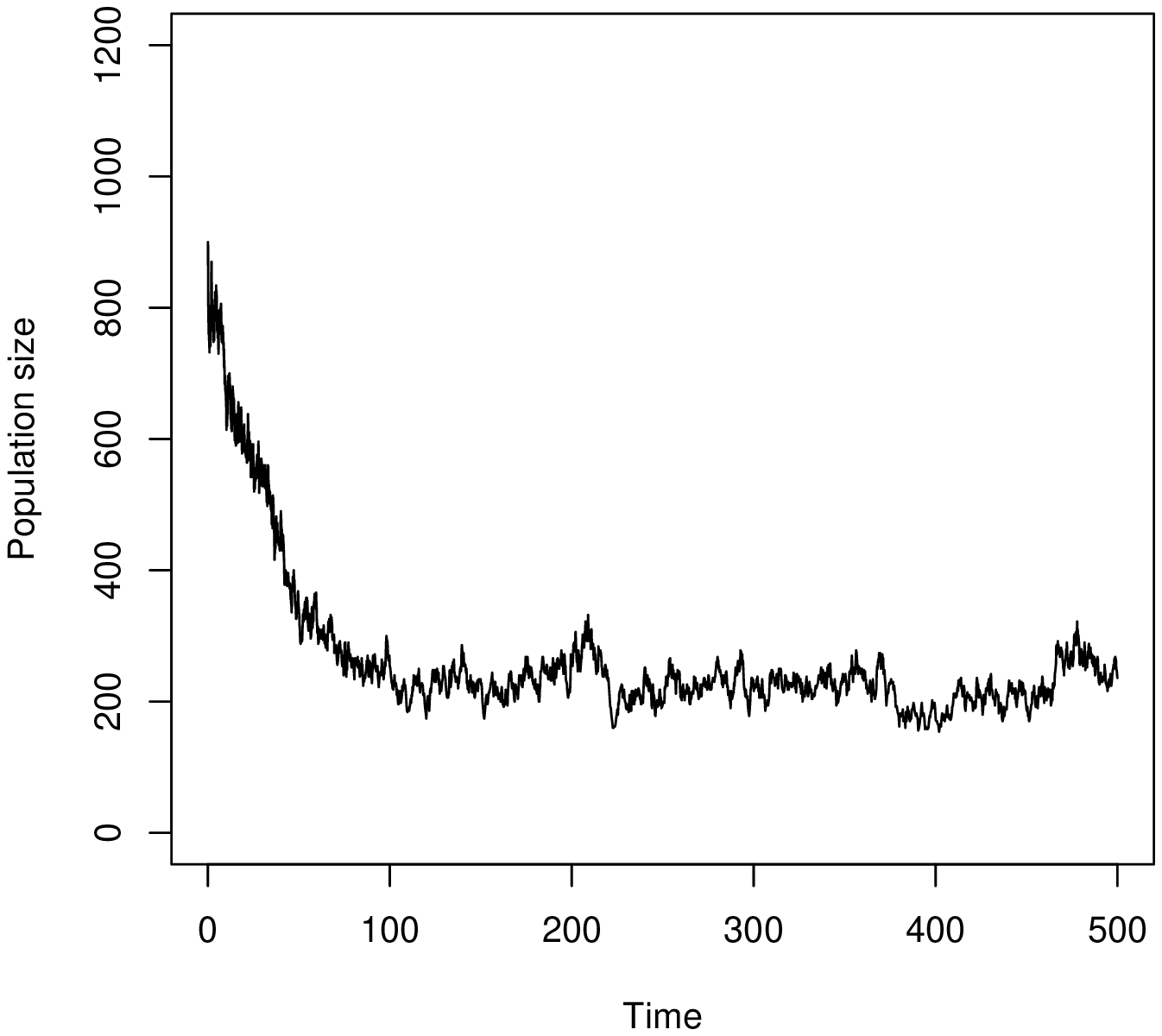}
\end{minipage}
\end{tabular}
\vspace{0.3cm}
\caption{{\small \textit{Population size dynamics as size at birth, growth rate and age at maturity coevolve. See Fig. \ref{figsim11} and main text for details. (a) $d_0=3,3.10^{-3}$, (b) $d_0=0.5$}}}
\label{figsize}
\end{center}
\end{figure}

\begin{figure}[!ht]
\vspace{1cm}
\begin{center}
\begin{tabular}[!ht]{cc}
\hspace{1.5cm}
\begin{minipage}[b]{.4\textwidth}\centering
\includegraphics[width=0.52\textwidth,height=0.24\textheight,angle=0,trim= 8cm 3cm 2cm 2cm]{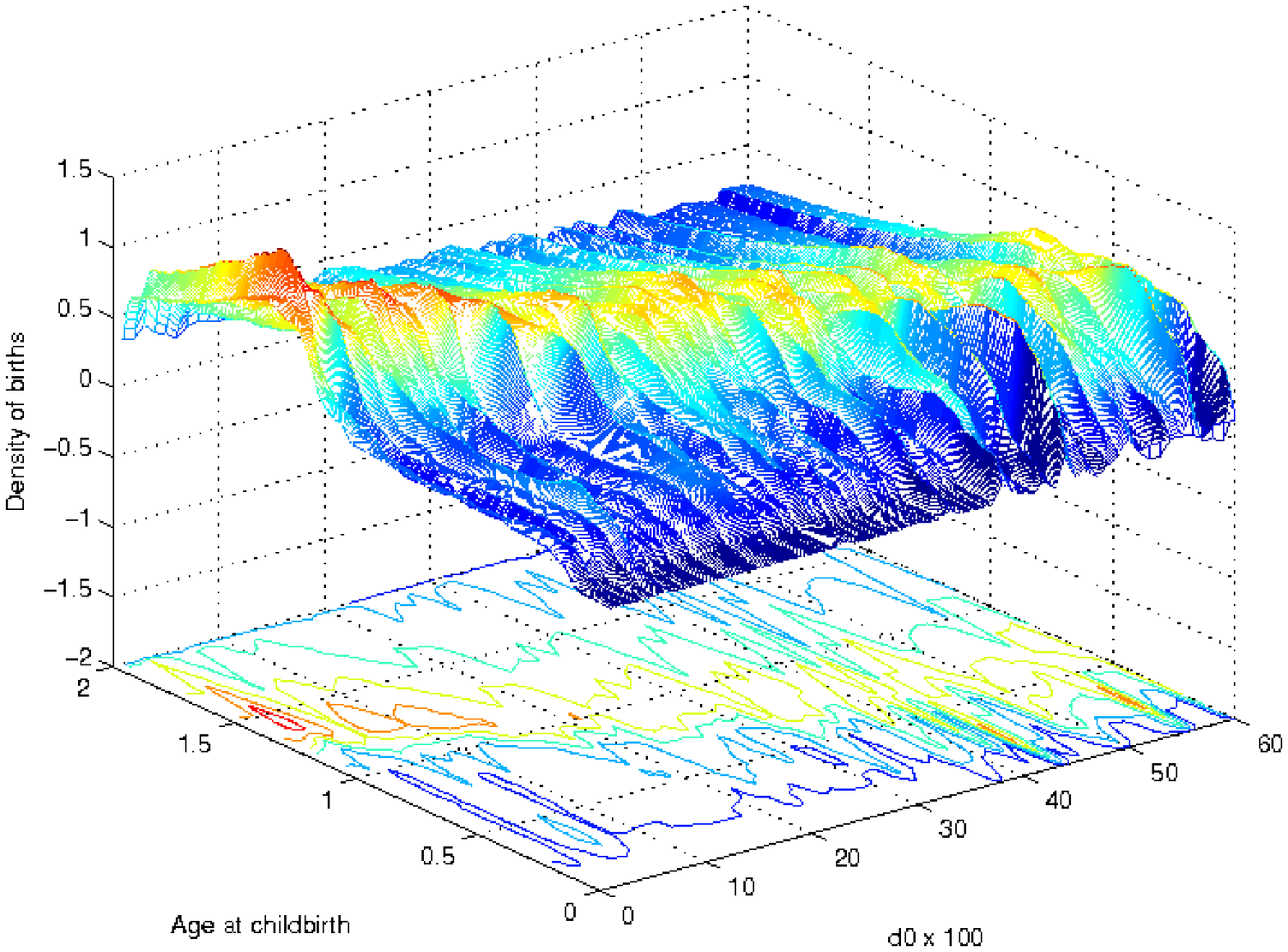}
\vspace{2cm}
\end{minipage}\hspace{1.4cm} &
\begin{minipage}[b]{.4\textwidth}\centering
\includegraphics[width=0.52\textwidth,height=0.24\textheight,angle=0,trim= 8cm 3cm 2cm 2cm]{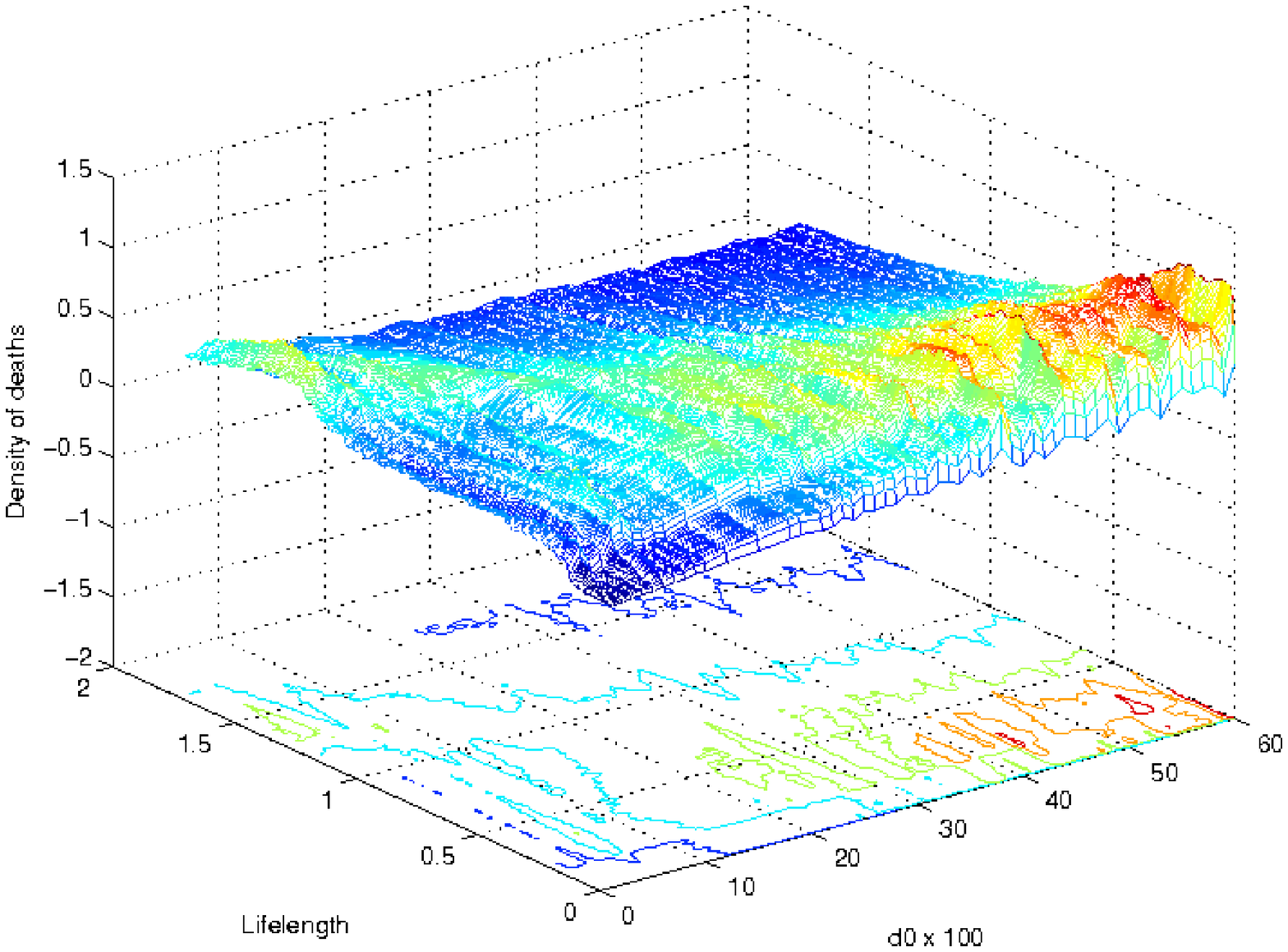}
\vspace{2cm}
\end{minipage}\\
(a) &
(b) \\
\end{tabular}
\vspace{-0.2cm}
\caption{{\small \textit{Long term density of maternal age at birth (a) and lifespan (b) with respect to the mortality cost of growth, $d_0$. Horizontal projection contours show that as $d_0$ increases, the distribution of maternal age at birth changes smoothly (a), whereas the distribution of lifespan undergoes an abrupt transition (around $d_0=0.2$) (b). See main text for details. }}}\label{figsim12}
\end{center}
\end{figure}

 When $d_0$ is large (Fig.\ref{figsim11}(b), $d_0=0.5$), a state with moderate growth and early maturation evolves. We observe that the growth rate decreases towards $1$ and evolution seeks to minimize the age at maturity ($a_M\simeq 0.3$). To understand the changes when $d_0$ increases, let us introduce the probability that an individual of traits $(x_0,g,a_M)$ born at $t=c$ survives until age at maturity $a_M$ in a population $(Z_t)_{t\in \R_+}$:
\begin{eqnarray}
\Pi(x_0,g,a_M)
 =  \mathbb{E}\left[\exp\left(-\int_{c}^{c+a_M} d\big(x_0,g,a_M,a,Z_{a} \big)da\right)\right].\label{probasurvie}
\end{eqnarray}
When the population size is bounded by $N$, the probability $\Pi(x_0,g,a_M)$ is upper and lower bounded by:
\begin{equation}\exp\left(-a_M\left(d_0g + N \cdot 10^{-5}\frac{2}{300}\right)\right)\leq \Pi(x_0,g,a_M)\leq \exp\left(-d_0 g a_M\right).\label{encadrementsurvie}\end{equation}In the simulations, we observe that the size of the population remains bounded and that it stabilizes around an "equilibrium" value (600 for $d_0=3.3\,10^{-3}$ and 200 for $d_0=0.5$, see Fig. \ref{figsize}). The "equilibrium" traits obtained in the two simulations of Figure \ref{figsim11} are distributed around $(x_0,g,a_M)=(1.8,2,1.3)$ and $(x_0,g,a_M)=(2.6,0.8,0.3)$. The strategy that is observed can be naturally explained. For $d_0=0.5$ and $N=600$, we hence obtain $0.27251\leq \Pi(1.8, 2,1.3)\leq 0.27254$, instead of $ \Pi(1.8, 2,1.3)\simeq 0.99$ when $d_0=3,3.10^{-3}$, whereas $0.88690 \leq \Pi(2.6,0.8,0.3)\leq 0.88692$. For $d_0=0.5$, the traits $(x_0,g,a_M)=(2.6,0.8,0.3)$ are more competitive. They correspond to individual who have a higher probability to reach maturity and reproduce, even if these individuals are likely to be exposed to higher competition pressure that individuals with trait $(x_0,g,a_M)=(1.8, 2,1.3)$.

\par In the case of a favorable environment with abundant ressources ($d_0$ is small), the growth period has a small cost and a long growth phase ($a_M$ large) allows individuals to reach a \textit{mortality refuge}: \ie a state where they escape the effects of strong competition pressure. When $d_0$ increases, the mortality refuge becomes inaccessible: long period of initial growth result in low probabilities of survival until reproduction. An alternative life profile appears, with early reproduction (small $a_M$) and small individuals. Similar phenomena have been discovered and discussed in Taborsky \textit{et al.} \cite{taborskydieckmannheino}.

\par The transition from "cheap" growth (low $d_0$, Fig. \ref{figsim11}(a)) to "expensive" growth (high $d_0$, Fig. \ref{figsim11}(b)) is investigated numerically in Fig. \ref{figsim12}, where $d_0$ varies between 0.01 and 0.6. Figure \ref{figsim12}(a) shows the kernel density estimator (\eg \cite{bosqlecoutre}) of the maternal age at birth. Fig. \ref{figsim12}(b) shows the density estimator of lifespan. When $d_0$ increases from $0.01$ (favorable environment) to $0.6$ (hostile environment), lifespan and reproduction age decrease continuously from 1.5 to 0.5. Compared to the results of \cite{taborskydieckmannheino}, the surfaces displayed in Fig. \ref{figsim12} do not show strict discontinuities, but the smooth variation of the distribution of maternal age at birth (a) contrasts with an abundant shift in the distribution of lifespan (b). This shift occurs in spite of the survival probability (\ref{probasurvie}) (and hence the lifespan density) being a continuous function of $d_0$.

\section{Conclusion}

\par Classical models for the dynamics of populations with continuous age structure ignore stochastic processes operating at the level of individuals, and therefore take the form of deterministic PDEs, \eg the celebrated McKendrick-von Foerster equation. Stochastic models that start at the level of individual processes and account for interactions between individuals face serious mathematical challenges. One way forward is to derive large population limits from these IBMs. This approach was followed here for a large class of IBMs, for which the large population approximation is a PDE that generalizes the Mc Kendrick-Von Foerster model. There are qualitative limitations of the large population approximation, that we demonstrated by showing almost sure extinction of the logistic age-structured population process. However, stochastic and deterministic approaches appear to yield complementary insights into the population dynamics: prior to extinction, stochastic trajectories are proved to spend an exponentially distributed time near the solution of the deterministic approximation model. In general, the rigorous derivation of the deterministic limit from stochastic individual-level processes allows one to derive confidence intervals for model parameters and thus open new ways of confronting models to real individual data.

\par In populations that are structured by age and trait variation (due \eg to genetic mutation), the analytical study even of the large population approximation becomes intractable. But the rigorous construction of the individual-based stochastic process yields an efficient algorithm for numerical simulations of population age and trait distributions. Two biological examples were presented. In the first example, offspring size varies genetically and evolves under size-dependent competition and the genetic constraint of a tradeoff with the birth rate. Individual growth and the resulting age and size structures have dramatic influences on trait evolution. Moderate growth shapes the population size structure in a way that exacerbates competition and favors the rapid split of the population into two 'evolutionary branches' (\ie the trait distribution becomes bimodal). With more rapid individual growth, the size distribution widens even more and competition intensifies to the point where the divergence of trait branches becomes impossible (\ie the population trait distribution remains unimodal).
\par In the second example, three traits were allowed to vary genetically and co-evolve: offspring size, growth rate, and age at maturity. The corresponding population model assumes size- and stage-dependent competition. As in the previous example, offspring size evolves under the genetic constraint of a trade-off with the birth rate. The growth rate evolves under a risk-competition tradeoff : faster growth (requiring \eg acquisition of more resources, hence riskier behavior) entails a mortality cost but ensures a larger size at maturity and hence a competitive advantage in the reproductive stage. Age at maturity evolves under a reproduction-competition tradeoff: reproductive (\ie mature) individuals face a higher mortality risk due to competition than juveniles. Model simulations were performed to analyze the effect of the mortality cost of growth on the traits' coevolution. High cost promotes the evolution of a moderate growth rate, together with large size at birth and very early maturity. The corresponding distribution of longevity in the population is skewed towards low values. Interestingly, as the mortality cost of growth decreases, the longevity distribution shows a relatively abrupt shift towards large values - a qualitative change in the population demography that was not anticipated given the continuous dependence of the longevity density on the cost of growth.
\par These examples highlight that the mathematical and numerical analysis of stochastic population models with age and trait structure have the potential to uncover unexpected phenomena of biological interest and to advance ecological and evolutionary population theory in significant ways.

{\small
\bibliographystyle{plain}
\bibliography{biblio}
}
\end{document}